                         \def\version{March 6, 2002}                         %

\documentclass[reqno,11pt]{amsart}

\usepackage{amsmath,epsfig, amsfonts}
\usepackage{amssymb}

\newfam\Bbbfam
\font\tenBbb=msbm10 \font\sevenBbb=msbm7 \font\fiveBbb=msbm5
\textfont\Bbbfam=\tenBbb \scriptfont\Bbbfam=\sevenBbb
\scriptscriptfont\Bbbfam=\fiveBbb

\renewcommand{\r}     {{{r}}}    

\newcommand{\ev}    {{\Ecal}}

\renewcommand{\i}   {{\operatorname{i}}}

\newcommand{\R}     {{\mathbb{R}}}
\newcommand{\Z}     {{\mathbb{Z}}}
\newcommand{\N}     {{\mathbb{N}}}

\newcommand{\Ssig}    {S}

\def\1{{\mathchoice {1\mskip-4mu\mathrm l}
{1\mskip-4mu\mathrm l} {1\mskip-4.5mu\mathrm l}
{1\mskip-5mu\mathrm l}}}

\def\comment#1{}
\newtheoremstyle{thm}{2ex}{2ex}{\itshape\rmfamily}{}
{\bfseries\rmfamily}{}{1.7ex}{}

\newtheoremstyle{rem}{1.3ex}{1.3ex}{\rmfamily}{}
{\itshape\rmfamily}{}{1.5ex}{}

\newenvironment{proofsect}[1]
{\vskip0.1cm\noindent{\bf #1}}{\vspace{0.15cm}}

\newtheorem{theorem}{Theorem}[section]
\newtheorem{lemma}[theorem]{Lemma}
\newtheorem{prop}[theorem] {Proposition}
\newtheorem{cor}[theorem]  {Corollary}

\newtheorem{step}{STEP}

\newcommand{\en}       {\end{equation}}
\newcommand{\eq}       {\begin{equation}}

\newcommand{\eqry}     {\begin{eqnarray}}
\newcommand{\enqry}    {\end{eqnarray}}
\newcommand{\eqarray}  {\begin{eqnarray}}
\newcommand{\enarray}  {\end{eqnarray}}
\newcommand{\eqarraystar} {\begin{eqnarray*}}
\newcommand{\enarraystar} {\end{eqnarray*}}
\newcommand{\bel}      {\begin{lemma}}
\newcommand{\el}       {\end{lemma}}

\newcommand{\bec}      {\begin{cor}}
\newcommand{\ec}       {\end{cor}}

\newcommand{\bes}      {\begin{step}}
\newcommand{\es}       {\end{step}}
\newcommand{\bea}      {\begin{array}}
\newcommand{\ea}       {\end{array}}
\newcommand{\bpr}      {\begin{proof}}
\newcommand{\epr}      {\end{proof}}

\newcommand{\da}       {\downarrow}

\renewcommand{\section}{\secdef\sct\sect}
\newcommand{\sct}[2][default]{\refstepcounter{section}
\vspace{0.5cm} \setcounter{equation}{0}
\centerline{\scshape \arabic{section}.\ #1} \vspace{0.3cm}}
\newcommand{\sect}[1]{\vspace{0.5cm}
\centerline{\large\scshape #1} \vspace{0.3cm}}

\renewcommand{\subsection}{\secdef\subsct\sbsect}
\newcommand{\subsct}[2][default]{\refstepcounter{subsection}
\nopagebreak \vspace{0.5\baselineskip} {\flushleft\bf
\arabic{section}.\arabic{subsection}~\bf #1 } \nopagebreak}
\newcommand{\sbsect}[1]{\vspace{0.1cm}\noindent
{\bf #1}\vspace{0.1cm}}

\renewcommand{\subsubsection}{\secdef\subsubsect\sbsbsect}
\newcommand{\subsubsect}[2][default]{\refstepcounter{subsubsection}
\nopagebreak \vspace{0.1\baselineskip} {\flushleft\sffamily\slshape
\arabic{section}.\arabic{subsection}.\arabic{subsubsection}
\sffamily #1\/.} }
\newcommand{\sbsbsect}[1]{\vspace{0.1cm}\noindent
{\bf #1}\ }

\newcommand{\eps}    {\varepsilon}

\renewcommand{\d}      {{\rm d}}

\newcommand{\heap}[2]{\genfrac{}{}{0pt}{}{#1}{#2}}

\newcommand{\Dcal}   {{\mathcal D }}
\newcommand{\Ecal}   {{\mathcal E }}

\newcommand{\Ocal}   {{\mathcal O }}

\newcommand{\smallsup}[1] {{\scriptscriptstyle{({#1}})}}

\begin{document}

\title[One-dimensional weakly interacting polymers]
{\large Weak interaction limits\\
for one-dimensional random polymers}

\author[Remco van der Hofstad, Frank den Hollander, Wolfgang K\"onig]{}
\maketitle
\thispagestyle{empty}
\vspace{1cm}

\centerline{\small \version}

\vspace{0.5cm}

\begin{center}
Remco van der Hofstad
\footnote{Department of Applied Mathematics,
Delft University of Technology, Mekelweg 4,
2628 CD Delft, The Netherlands.}
\footnote
{Present address: Department of Mathematics and Computer Science,
Eindhoven University of Technology, P.O.\ Box  513,
5600 MB Eindhoven, The Netherlands.
{\tt rhofstad@win.tue.nl}}\\
Frank den Hollander
\footnote{EURANDOM, P.O.\ Box 513,
5600 MB Eindhoven, The Netherlands.
{\tt denhollander@eurandom.tue.nl}}\\
Wolfgang K\"onig
\footnote{Institut f\"ur Mathematik, TU Berlin,
Stra{\ss}e des 17.\ Juni 136, D-10623 Berlin, Germany.
{\tt koenig@math.tu-berlin.de}}
\end{center}

\vspace{2cm}

\begin{quote}
{\small
{\bf Abstract:}
In this paper we present a new and flexible method to show that, in one dimension,
various self-repellent random walks converge to self-repellent Brownian motion in
the limit of weak interaction after appropriate space-time scaling. Our method is
based on cutting the path into pieces of an appropriately scaled length, controlling
the interaction between the different pieces, and applying an invariance principle
to the single pieces. In this way we show that the self-repellent random walk large
deviation rate function for the empirical drift of the path converges to the
self-repellent Brownian motion large deviation rate function after appropriate
scaling with the interaction parameters. The method is considerably simpler than
the approach followed in our earlier work, which was based on functional analytic
arguments applied to variational representations and only worked in a very limited
number of situations.

We consider two examples of a weak interaction limit: (1) vanishing
self-repellence, (2) diverging step variance. In example (1), we recover
our earlier scaling results for simple random walk with vanishing
self-repellence and show how these can be extended to random walk with
steps that have zero mean and a finite exponential moment. Moreover, we
show that these scaling results are stable against adding self-attraction,
provided the self-repellence dominates. In example (2), we prove a conjecture
by Aldous for the scaling of self-avoiding walk with diverging step variance.
Moreover, we consider self-avoiding walk on a two-dimensional horizontal
strip such that the steps in the vertical direction are uniform over the
width of the strip and find the scaling as the width tends to infinity.
}
\end{quote}

\vfill

\bigskip\noindent
{\it 2000 Mathematics Subject Classification.} 60F05, 60F10, 60J55, 82D60.

\medskip\noindent
{\it Keywords and phrases.} Self-repellent random walk and Brownian motion,
invariance principles, large deviations, scaling limits, universality. \eject

\setcounter{section}{0}


\section{Polymer measures}
\label{sec-WSAW}

A polymer is a long chain of atoms or molecules, often referred to as monomers,
which have a tendency to repel each other. This self-repellence comes from the
excluded-volume-effect: two molecules cannot occupy the same space. The
self-repellence causes the polymer to spread itself out more than it would do
in the absence of self-repellence. The most widely used ways to describe a
polymer are the {\it Domb-Joyce model\/}, respectively, the {\it Edwards model\/},
which start from random walk, respectively, Brownian motion and build in an
appropriate penalty for self-intersections. In Sections~\ref{model} and
\ref{sec-EM} we introduce these two models (in dimension one) and list some
known results about their space-time scaling. In Section~\ref{results} we
consider a number of variations on the Domb-Joyce model and formulate our main
results, which are weak interaction limits showing that all these models
scale to the Edwards model in the limit of weak interaction. Section \ref{sec-LD}
reviews some large deviation results for the Domb-Joyce model and the Edwards model,
while Sections~\ref{sec-proof(i)}--\ref{sec-proofs} contain the proofs of the
theorems in Section~\ref{results}. In Section~\ref{sec-disc} we close with a
brief discussion of the method of proof and of some open ends.

A general background on polymers from a physics and chemistry point of view
may be found in \cite{vdZ}, a survey of mathematical results for one-dimensional
polymers appears in \cite{HK01}.

\subsection{The Domb-Joyce model.}
\label{model}

Let $(S_n)_{n\in\N_0}$ be a random walk on $\Z$ starting at the
origin $(S_0=0)$. Let $P$ be the law of this random walk and let
$E$ be expectation with respect to $P$. Assume that the random walk is
irreducible and that
    \begin{equation}
    \label{mv}
    E(S_1) = 0, \qquad E(e^{\eps |S_1|})<\infty
    \quad \mbox{for some }\eps>0.
    \end{equation}
Throughout the paper,
    \begin{equation}
    \label{sigma}
    \sigma^2 = E|S_1|^2 \in (0,\infty)
    \end{equation}
denotes the step variance.

Fix $n\in\N$, introduce a parameter $\beta\in [0,\infty]$, and
define a probability law $Q_n^{\beta}$ on $n$-step paths by setting
    \begin{equation}
    \label{PM}
    \frac{\d Q_n^{\beta}}{\d P}[\cdot]
    =\frac{1}{Z_n^{\beta}}e^{-\beta H_n[\cdot]},\qquad
    Z_n^{\beta}=E(e^{-\beta H_n}),
    \end{equation}
with
    \begin{equation}
    \label{Hamilt}
    H_n\left[(S_i)_{i=0}^n\right] =  \sum_{\heap{i,j=0}{i\not= j}}^n
    \1_{\{S_i=S_j\}}=\sum_{x\in\Z}\ell_n(x)^2 - (n+1)
    \end{equation}
the intersection local time up to time $n$, where
    \begin{equation}
    \label{deflt}
    \ell_n(x) = \#\{0\leq i \leq n\colon S_i=x\}, \qquad x\in \Z,\
    \end{equation}
is the local time at site $x$ up to time $n$. The law $Q_n^{\beta}$ is called
the $n$-polymer measure with strength of self-repellence $\beta$. The path
receives a penalty $e^{-2\beta}$ for every self-intersection. The term
$n+1$ in \eqref{Hamilt} can be trivially absorbed into the normalization.

In the case $\beta=\infty$, with the convention $e^{-\infty H_n}=\1_{\{H_n=0\}}$,
the path measure $Q_n^{\infty}$ is the conditional probability law given that
there are no self-intersections up to time $n$, i.e., $Q_n^\infty=P(\,\cdot\mid H_n=0)$.
If single steps are equally probable under $P$, then $Q_n^\infty$ is the uniform
distribution on all $n$-step self-avoiding paths having a strictly positive probability
under $P$. The law $Q_n^\infty$ is known as the {\em self-avoiding walk}, and is
trivial for simple random walk but non-trivial when the random walk can make
larger steps.
\footnote{For $\beta\in(0,\infty)$, $Q_n^\beta$ is sometimes referred to as the
{\em weakly self-avoiding walk}.}

For the special case where
    \begin{equation}
    S_1 \mbox{ is symmetric with support } \{-L,\dots,-1,1,\dots,L\}
    \mbox{ for some } L \in \N,
    \label{rcond}
    \end{equation}
the following is known.

    \medskip
    \begin{theorem}[CLT and partition function]
    \label{thm-DJ*}
    Fix $\beta \in [0,\infty]$, assume~\eqref{rcond}, and exclude the trivial
    case $(\beta,L)=(\infty,1)$. Then there are numbers $r^*,\theta^*,\sigma^*
    \in (0,\infty)$ (depending on $\beta$ and on the distribution of $S_1$)
    such that:
    \begin{enumerate}
    \item[(i)]
    Under the law $Q_n^{\beta}$, the distribution of the scaled
    and normalized endpoint $(|S_n|-\theta^* n)/\sigma^*\sqrt n$ converges
    weakly to the standard normal distribution.
    \item[(ii)]
    $\lim_{n\to\infty} \frac 1n\log Z_n^{\beta}=-\r^*$.
    \end{enumerate}
    \end{theorem}

\medskip\noindent
Theorem~\ref{thm-DJ*}(i) is contained in \cite[Theorem 1.1]{K3}, Theorem~\ref{thm-DJ*}(ii)
is proved in \cite{K2} for $\beta<\infty$ and in \cite{K1} for $\beta=\infty$. For
$L=1$, the law of large numbers contained in Theorem~\ref{thm-DJ*}(i) first appeared
in Greven and den Hollander \cite{GH}.

\subsection{The Edwards model.}
\label{sec-EM}

Let $B=(B_t)_{t \geq 0}$ be a standard Brownian motion on $\R$ starting at the origin
($B_0=0$). Let $\widehat P$ be the Wiener measure and let $\widehat E$ be expectation
with respect to $\widehat P$. For $T>0$ and $\beta\in[0,\infty)$, define a probability
law $\widehat Q_T^{\beta}$ on paths of length $T$ by setting
    \begin{equation}
    \label{Ed mod}
    \frac{\d\widehat Q_T^\beta}{\d\widehat P}[\cdot]
    =\frac{1}{\widehat Z_T^\beta}e^{-\beta\widehat H_T[\cdot]},
    \qquad \widehat Z_T^\beta=\widehat E(e^{-\beta\widehat H_T}),
    \end{equation}
with
    \begin{equation}
    \label{HT}
    \widehat H_T\left[(B_t)_{t\in[0,T]}\right]
    =\int_0^T{\rm d}u \int_0^T {\rm d}v\,\,\delta(B_u-B_v)
    =\int_{\R}L(T,x)^2\,{\rm d}x
    \end{equation}
the Brownian self-intersection local time up to time $T$. The middle
expression in (\ref{HT}) is formal only. In the last expression
the Brownian local times $L(T,x)$, $x\in\R$, appear. The law $\widehat Q_T^\beta$
is called the $T$-polymer measure with strength of self-repellence $\beta$.
The Brownian scaling property implies that
    \begin{equation}
    \bigl(L(t,x)\bigr)_{t\in[0,T],x\in\R} \stackrel{\Dcal}{=}
    \left(\beta^{-\frac 13}L(\beta^{\frac 23}t,\beta^{\frac 13}x)
    \right)_{t\in[0,T],x\in\R}, \qquad \beta,T>0
    \end{equation}
(here $\stackrel{\Dcal}{=}$ means equal in distribution under
$\widehat P$), and hence that
    \begin{equation}
    \label{scalEM}
    \widehat Q_T^\beta\left((B_t)_{t\in[0,T]}\in\cdot~\right)
    = \widehat Q_{\beta^{\frac23} T}^1\left((\beta^{-\frac13}
    B_{\beta^{\frac 23}t})_{t\in[0,T]}\in \cdot~\right), \qquad \beta,T>0.
    \end{equation}

\medskip
\begin{theorem}[CLT and partition function]
\label{thm-CLTEDW} There are numbers $a^*, b^*, c^*\in(0,\infty)$ such that,
for any $\beta\in(0,\infty)$:
\begin{enumerate}
\item[(i)]
    Under the law $\widehat Q_T^\beta$, the distribution of the scaled and
    normalized endpoint $(|B_T|-b^*\beta^{\frac 13} T)/c^*\sqrt T$ converges
    weakly to the standard normal distribution.
\item[(ii)]
    $\lim_{T\to\infty}\frac 1T\log \widehat Z_T^\beta=-a^*\beta^{\frac 23}$.
\end{enumerate}
\end{theorem}

\medskip\noindent
Theorem~\ref{thm-CLTEDW} is proved in \cite{HHK1}. Rigorous bounds on $a^*$,
$b^*$, $c^*$ appeared in \cite[Theorem 3]{Hconst}. The numerical values are:
$a^*\approx 2.19$, $b^*\approx 1.11$, $c^*\approx 0.63$. The law of large numbers
contained in Theorem~\ref{thm-CLTEDW}(i) first appeared in Westwater \cite{We3}.

\newpage


\section{Main results}
\label{results}

In this section we formulate and explain our main results, all of which are
weak interaction limits for the large space-time scaling of the one-dimensional
Domb-Joyce model introduced in Section~\ref{model} and various related models.
In all cases the scaling is the same as that of the Edwards model introduced
in Section~\ref{sec-EM}, showing that universality holds. Two examples of a
weak interaction limit are considered: $\beta\da 0$ and $\sigma\to\infty$.

Section~\ref{mainres} considers the Domb-Joyce model, Section~\ref{sec-attr}
the Domb-Joyce model with added self-attraction, and Section~\ref{twodim}
self-avoiding walk on a two-dimensional strip. In Section~\ref{weaklim}
we describe some invariance principles that are needed in the proofs appearing
in Sections~\ref{sec-proof(i)}--\ref{sec-proofs}. A brief discussion of our results
and our method of proof can be found in Section~\ref{sec-disc}.

\subsection{Two weak interaction limits for self-repellent polymers.}
\label{mainres}

Consider an arbitrary random walk $(\Ssig_n)_{n\in\N_0}$ on $\Z$ satisfying \eqref{mv},
respectively, the two technical conditions (\ref{sigmares}--\ref{expmomsigma})
introduced in Section~\ref{weaklim}.

    \medskip
    \begin{theorem}[LLN]
    \label{ratescal}
    $\mbox{}$
    \begin{enumerate}
    \item[(i)]
    Fix $\sigma\in(0,\infty)$. Then, under \eqref{mv},
    \begin{equation}
    \label{LLN}
    \lim_{\beta\downarrow 0}\limsup_{n\rightarrow \infty}
    Q_n^{\beta}\Bigl(\Big|\frac{|\Ssig_n|}{\beta^{\frac 13}n}
    -b^*\sigma^{\frac 23}\Big| \geq \eps\Bigl)=0\qquad \forall \eps>0.
    \end{equation}
    \item[(ii)]
    Fix $\beta=\infty$. Then, under \eqref{sigmares}--\eqref{expmomsigma},
    \begin{equation}
    \lim_{\sigma\rightarrow \infty}\limsup_{n\rightarrow \infty}
    Q_n^{\infty}\Bigl(\Big|\frac{|\Ssig_n|}{\sigma^{\frac 23}n}
    -b^*\Big| \geq \eps\Bigl)=0\qquad \forall \eps>0.
    \end{equation}
    \end{enumerate}
    \end{theorem}

\medskip\noindent
Theorem~\ref{ratescal} is proved in Sections~\ref{sec-proof(i)}--\ref{sec-proof(ii)}.
It is to be viewed as an {\it approximative law of large numbers\/} for the endpoint
$\Ssig_n$ of the polymer, since it states that the asymptotics of $|\Ssig_n|/n$ as
$n\to\infty$ behaves like $b^*\sigma^{\frac 23}\beta^{\frac 13}$ as $\beta\da 0$,
respectively, like $b^*\sigma^{\frac 23}$ as $\sigma\to\infty$. Note that in
Theorem~\ref{ratescal}(i) the asymptotics does not depend on the details of the
random walk other than its step variance.

In the special case of \eqref{rcond}, where the central limit
theorem is known (recall Theorem \ref{thm-DJ*}(i)), we obtain the
following two corollaries for the scaling of the parameters $r^*$
and $\theta^*$ as $\beta\downarrow 0$, respectively,
$\sigma\to\infty$. To stress this dependence, we write
$r^*=r^*(\beta), \theta^*=\theta^*(\beta)$. Both these corollaries
are also proved in
Sections~\ref{sec-proof(i)}--\ref{sec-proof(ii)}.

    \medskip
    \begin{cor}[Scaling rate and drift]
    \label{asy}
    Fix $\sigma\in(0,\infty)$. Then, under \eqref{rcond},
    \begin{equation}
    \label{asyvalue}
    \r^*(\beta)\sim a^*\sigma^{-\frac 23}\beta^{\frac 23}, \qquad
    \theta^*(\beta)\sim b^*\sigma^{\frac 23}\beta^{\frac 13},
    \qquad \beta\da 0.
    \end{equation}
    \end{cor}

\medskip\noindent
For the nearest-neighbor random walk ($\sigma^2=1$), the assertions in
Corollary~\ref{asy} were already proved in \cite[Theorems 4--6]{HH}.
However, the proof used heavy functional analytic
tools and gave no probabilistic insight. For $\sigma^2>1$ this route seems
inaccessible, so it is nice that here the scaling comes out more generally.

    \medskip
    \begin{cor}[Scaling rate and drift]
    \label{asy*}
    Fix $\beta=\infty$. Then, under \eqref{rcond} and
    \eqref{sigmares}--\eqref{expmomsigma},
    \begin{equation}
    \r^*(\infty)\sim a^*\sigma^{-\frac 23}, \qquad
    \theta^*(\infty)\sim b^*\sigma^{\frac 23},
    \qquad \sigma\to\infty.
    \end{equation}
    \end{cor}

\medskip\noindent
The second assertion in Corollary~\ref{asy*} settles a conjecture due to Aldous
\cite[Section 7(B)]{Al}, although Aldous misses the factor $b^*$.

We believe that also
\begin{equation}
\label{sigscal}
\sigma^*(\beta) \to c^*, \quad \beta\da 0, \qquad \mbox{respectively}
\qquad \sigma^*(\infty) \to c^*, \quad \sigma\to\infty,
\end{equation}
but we are unable to prove this. The reason why will become clear in
Section~\ref{sec-proofcor1}. For nearest-neighbor random walk, the first assertion in
\eqref{sigscal} was proved in \cite{HHK2}.

Our approach is flexible enough to allow for a coupled limit $n\to\infty$ and
$\beta\da 0$, respectively, $\sigma\to\infty$.

    \medskip
    \begin{theorem}[Coupled LLN]
    \label{thm-betan}
    $\mbox{}$
    \begin{itemize}
    \item[(i)]
    Fix $\sigma\in(0,\infty)$, and assume \eqref{mv}. If $\beta$ is
    replaced by $\beta_n$ satisfying $\beta_n\to 0$ and $\beta_n n^{\frac 32}
    \to\infty$ as $n\to\infty$, then
    \begin{equation}
    \lim_{n\rightarrow \infty}
    Q_n^{\beta_n}\Bigl(\Big|\frac{|\Ssig_n|}{\beta_n^{\frac 13}n}
    -b^*\sigma^{\frac 23}\Big|\geq\eps\Bigl)=0 \qquad \forall\eps>0.
    \end{equation}
    \item[(ii)]
    Fix $\beta=\infty$, and assume \eqref{rcond} and \eqref{sigmares}--\eqref{expmomsigma}.
    If $\sigma$ is replaced by $\sigma_n$ satisfying $\sigma_n\to\infty$ and
    $\sigma_n n^{-\frac 32}\to 0$ as $n\to\infty$, then
    \begin{equation}
    \lim_{n\rightarrow \infty}
    Q_n^{\infty}\Bigl(\Big|\frac{|\Ssig_n|}{\sigma_n^{\frac 23}n}
    -b^*\beta^{\frac 13}\Big|\geq\eps\Bigl)=0 \qquad \forall\eps>0.
    \end{equation}
    \end{itemize}
    \end{theorem}

\medskip\noindent
Theorem~\ref{thm-betan} is proved in Section~\ref{proof-betan}. For simple random
walk ($\sigma^2=1$), the assertion in Theorem~\ref{thm-betan}(i) was already proved in
\cite[Theorem 1.5]{HHK2}. Note that the conditions on $\beta_n$, respectively,
$\sigma_n$ keep the scaling out of the central limit regime.

\subsection{Weak interaction limit for self-repellent and self-attractive polymers.}
\label{sec-attr}

The method introduced in this paper extends to the situation where self-attraction
is added to the polymer. In \eqref{PM}, we replace $\beta H_n$ by
    \begin{eqnarray}
    \label{Hamiattr}
    H_n^{\beta,\gamma}&=&\beta\sum\limits_{\heap{i,j=0}{i\not= j}}^n
    \1_{\{S_i=S_j\}}-\frac{\gamma}{2}\sum\limits_{\heap{i,j=0}{i\not= j}}^n
    \1_{\{|S_i-S_j|=1\}}\nonumber\\
    &=&(\beta-\gamma)\sum\limits_{x\in\Z}\ell^2_n(x)
    +\frac \gamma {2}\sum\limits_{x\in\Z}[\ell_n(x)-\ell_n(x+1)]^2
    -\beta(n+1),
    \end{eqnarray}
where $\beta,\gamma\in(0,\infty)$ are parameters, and $(S_n)_{n\in\N_0}$ is an
arbitrary random walk on $\Z$ satisfying \eqref{mv}. In words, $H_n^{\beta,\gamma}$
is equal to $\beta$ times twice the number of self-intersections up to time $n$
minus $\gamma$ times twice the number of self-contacts up to time $n$. The law
$Q_n^{\beta,\gamma}$ gives a penalty $e^{-2\beta}$ to every pair of monomers
at the same site and a reward $e^{\gamma}$ to every pair of monomers at
neighboring sites. The term $\beta (n+1)$ in \eqref{Hamiattr} can again be
trivially absorbed into the normalization.

The scaling behavior under $Q_n^{\beta,\gamma}$ was studied (in arbitary dimension)
in \cite{HK00}. It was shown that there is a phase transition at $\beta=\gamma$, namely,
the polymer collapses on a finite (random) number of sites when $\gamma>\beta$, while
it visits order $n$ sites when $\gamma<\beta$. Furthermore, in dimension one, a law
of large numbers and a central limit theorem for the endpoint $S_n$ under
$Q_n^{\beta,\gamma}$, analogous to Theorem~\ref{thm-DJ*}(i), were derived under
the restriction $0<\gamma<\beta-\frac 12 \log 2$.

We want to obtain the analogue of Theorem~\ref{ratescal}(i). In Theorem
\ref{thm-attr} below we abbreviate
    \begin{equation}
    \label{betagammarestr}
    \lim_{\beta,\gamma} \qquad \mbox{ for } \qquad \beta,\gamma \da 0 \mbox{ such that }
    0<\gamma<\beta\mbox{ and } \gamma(\beta-\gamma)^{-\frac 23}\to 0,
    \end{equation}
and likewise for $\liminf$ and $\limsup$.

    \medskip
    \begin{theorem}[LLN]
    \label{thm-attr}
    Fix $\sigma\in(0,\infty)$. Then, under \eqref{mv},
    \begin{equation}
    \lim\limits_{\beta,\gamma}\limsup_{n\rightarrow \infty}
    Q_n^{\beta,\gamma}\Bigl(\Big|\frac{|\Ssig_n|}{(\beta-\gamma)^{\frac 13}n}
    -b^*\sigma^{\frac 23}\Big| \geq \eps\Bigl)=0\qquad \forall \eps>0.
    \end{equation}
    \end{theorem}

\medskip\noindent
Theorem~\ref{thm-attr} is proved in Section~\ref{proof-attr}. Note that no
law of large numbers is known for small $\beta,\gamma$. If
\begin{equation}
\theta^*(\beta,\gamma)=\lim_{n\to\infty}E_{Q_n^{\beta,\gamma}}
\Bigl(\frac{|\Ssig_n|}{n}\Bigr) \in (0,\infty)
\end{equation}
would exist for fixed $\beta,\gamma$, then we could deduce from Theorem~\ref{thm-attr}
that $\lim_{\beta,\gamma}(\beta-\gamma)^{-\frac 13}\theta^* (\beta,\gamma)
=b^*\sigma^{\frac 23}$.

We believe that Theorem~\ref{thm-attr} fails without the restrictions on $\beta,\gamma$
in (\ref{betagammarestr}). There is also a coupled limit version of Theorem~\ref{thm-attr}
analogous to Theorem~\ref{thm-betan}, but we refrain from writing this down.

\subsection{Weak interaction limit for self-avoiding polymers on a
two-dimensional strip.}
\label{twodim}

Let $(X_n)_{n\in\N_0}=(\Ssig_n,U^L_n)_{n\in\N_0}$ be a random walk on the
strip $\Z\times \{-L,\dots,L\}$, where $(\Ssig_n)_{n\in\N_0}$ is a random walk
on $\Z$ satisfying \eqref{mv}, and $(U^L_n)_{n\in\N_0}$ is an i.i.d.\ sequence,
independent of $(\Ssig_n)_{n\in\N_0}$, such that $U^L_0$ is uniformly distributed
on $\{-L,\dots,L\}$. For this two-dimensional random walk, define its self-avoiding
version by putting $Q_n^{\infty,L}(\cdot)=P^{L}(\cdot\mid H_n=0)$, where $P^{L}$
is the law of $(X_n)_{n\in\N_0}$ and
    \begin{equation}
    \label{taudef}
    H_n=\sum_{\heap{i,j=0}{i\not= j}}^n
    \1_{\{X_i=X_j\}}
    \end{equation}
is the intersection local time up to time $n$.

Theorem~\ref{thm-twodim} below identifies the asymptotics of the endpoint
of the {\it first\/} coordinate, $\Ssig_n$, under the law $Q_n^{\infty,L}$
in the limit as $n\to\infty$ followed by $L\to\infty$, and also when the two
limits are coupled.

    \medskip
    \begin{theorem}[LLN and coupled LLN]
    \label{thm-twodim}
    Fix $\sigma\in(0,\infty)$ and assume \eqref{mv}.
    \begin{enumerate}
    \item[(i)]
    Then
    \begin{equation}
    \lim\limits_{L\rightarrow \infty}\limsup_{n\rightarrow \infty}
    Q_n^{\infty,L}\Bigl(\Big|\frac{|\Ssig_n|}{(4L)^{-\frac 13}n}
    -b^*\sigma^{\frac 23}\Big| \geq \eps\Bigl)=0\qquad \forall \eps>0.
    \end{equation}
    \item[(ii)]
    If $L$ is replaced by $L_n$ satisfying $L_n\to \infty$ and
    $L_n n^{-\frac 32}\to 0$ as $n\to\infty$, then
    \begin{equation}
    \label{Lnscal}
    \lim_{n\rightarrow \infty}
    Q_n^{\infty,L_n}\Bigl(\Big|\frac{|\Ssig_n|}{(4L_n)^{-\frac 13}n}
    -b^*\sigma^{\frac 23}\Big| \geq \eps\Bigl)=0\qquad \forall \eps>0.
    \end{equation}
    \end{enumerate}
    \end{theorem}

\medskip\noindent
Theorem~\ref{thm-twodim} is proved in Section~\ref{proof-twodim}. In \cite{AJ90},
it is shown that
\begin{equation}
\theta^*(L)=\lim_{n\to\infty}E_{Q_n^{\infty,L}}\Bigl(\frac{|\Ssig_n|}{n}\Bigl)
\in (0,\infty)
\end{equation}
exists for fixed $L$. Therefore, we deduce from Theorem~\ref{thm-twodim}(i)
that $\lim_{L\to\infty}(4L)^{\frac 13}\theta^*(L)=b^*\sigma^{\frac 23}$.

We close this section by making a comparison with self-avoiding walk on $\Z^2$.
One of the prominent open problems for this process is the asymptotic analysis
of its endpoint. The conjecture is that the endpoint runs on scale $n^{\frac 34}$.
Now, interestingly, in Theorem~\ref{thm-twodim}(ii) it is precisely the choice
$L_n=n^{\frac 34}$ that makes the two coordinates $\Ssig_n$ and $U_n^{L_n}$ run
on the {\it same\/} scale $n^{\frac 34}$. This suggests that for $L_n=n^{\frac 34}$
the behavior on the strip is a reasonable qualitative approximation to the behavior
on $\Z^2$.

Let us try to make this argument a bit more precise by appealing to an adaptation
of the well-known Flory argument (see \cite[Section 2.2]{MS}). Let $S=(S_n)_{n\in\N_0}
=(S_n^{\smallsup{1}},S_n^{\smallsup{2}})_{n\in\N_0}$ be two-dimensional simple
random walk. We may assume that $S^{\smallsup{1}}=(S^{\smallsup{1}}_n)_{n\in\N_0}$
and $S^{\smallsup{2}}=(S^{\smallsup{2}}_n)_{n\in\N_0}$ are two independent
one-dimensional simple random walks.\footnote
{Indeed, the projections of $S^{\smallsup{1}}$ and $S^{\smallsup{2}}$
onto the lines with slope 1 and $-1$ in $\R^2$, respectively, are two independent
copies of one-dimensional simple random walk on $\sqrt{2}\,\Z$.}
We want to investigate the quantity
        \begin{eqnarray}
        Z_n^\infty(\nu)&=& P\Big(\bigcap_{\stackrel{i,j=0}{i\neq j}}^n
        \{S_i\not=S_j\}\cap\{|S_n|\asymp n^\nu\}\Big)
        \label{rewrite2}\nonumber\\
        &=& E^{\smallsup{1}}\Big(\1\{|S_n^{\smallsup{1}}|\asymp n^{\nu}\}
        P\Big(\bigcap_{\stackrel{i,j=0}{i\neq j}}^n
        \{S_i\not=S_j\}\cap\{|S_n^{\smallsup{2}}|\asymp n^\nu\}
        \Big|S^{\smallsup{1}}\Big)\Big),
        \end{eqnarray}
where $P$ is the law of $S$, $E^{\smallsup{1}}$ is expectation with respect
to $S^{\smallsup{1}}$, and $\nu>0$ is an exponent to be determined later. Denote
the local times of $S^{\smallsup{1}}$ by $\ell^{\smallsup{1}}_n(x)$, $x\in\Z$.
Note that $S^{\smallsup{1}}$ has $\ell^{\smallsup{1}}_n(x)[\ell^{\smallsup{1}}_n(x)-1]$
self-intersections at $x\in\Z$. In order that $S$ has no self-intersections,
$S^{\smallsup{2}}$ must avoid a self-intersection at the $\sum_{x\in\Z}
\ell^{\smallsup{1}}_n(x)[\ell^{\smallsup{1}}_n(x)-1]$ time pairs at which $S^{\smallsup{1}}$
has self-intersections. Now, let us make the {\it crude approximation\/} that
$S_i^{\smallsup{2}}$, $i=0,\dots,n$, are i.i.d.\ uniformly distributed on
$\{-|S^{\smallsup{2}}_n|,\dots,|S^{\smallsup{2}}_n|\}$. Then, on the event
$\{|S_n^{\smallsup{2}}|\asymp n^{\nu}\}$, the probability that a self-intersection
of $S^{\smallsup{2}}$ occurs at a given time pair $i\not= j$ at which $
S_i^{\smallsup{1}}=S_j^{\smallsup{1}}$ is $\asymp n^{-\nu}$. (The idea behind
the approximation is that for large $n$ most self-intersections occur when $|i-j|$
is large.) The resulting model is precisely the one investigated in Theorem
\ref{thm-twodim}(ii) with $L_n\asymp n^{\nu}$. For this choice, \eqref{Lnscal}
yields that $\{S_n^{\smallsup{2}}\asymp n^{1-\frac{\nu}{3}}\}$ is typical.
Putting $\nu=1-\frac 13 \nu$, we find $\nu=\frac 34$.

\subsection{Invariance principles and assumptions on variance scaling.}
\label{weaklim}

The proofs of our weak interaction limits in Sections~\ref{mainres}--\ref{twodim}
will be based on a number of invariance principles, which we describe now.
Let $(B^\sigma_t)_{t\geq 0}$ be a Brownian motion with generator $\frac 12 \sigma^2
\Delta$, and write $\widehat H^\sigma_T$ for its intersection local time and
$L^\sigma(T,x)$, $x\in\R$, for its local times up to time $T$.

\medskip\noindent
I.\ The first invariance principle we will rely on was put forward in
\cite[Theorem 1.3]{BSl}:
\footnote{In fact, \cite[Theorem 1.3]{BSl} applies only to simple random walk,
but an inspection of its proof reveals that it in fact holds in the generality
of our setting.}
    \begin{equation}
    \label{ntoinfty}
    \Bigl(n^{-\frac 12}(\Ssig_{\lfloor nt\rfloor})_{t\in[0,T]},n^{-\frac 32}
    H_{\lfloor nT\rfloor}\Bigr)\stackrel{n\to\infty}
    {\Longrightarrow}\bigl((B^\sigma_t)_{t\in[0,T]},\widehat H^\sigma_T\bigr),
    \qquad \sigma,T>0.
    \end{equation}
This says that the Domb-Joyce model (for the random walk with variance $\sigma^2$)
at time $nT$ with strength of self-repellence $\beta n^{-\frac 32}$ converges, after
appropriate space-time scaling, to the Edwards model (for the Brownian motion with
generator $\frac 12 \sigma^2\Delta$) at time $T$ with strength of self-repellence
$\beta$. Another version of the same invariance principle is the assertion
    \begin{equation}
    \label{betatozero}
    \Bigl(\beta^{\frac 13}(\Ssig_{\lfloor\beta^{-\frac 23}t\rfloor})_{t\in[0,T]},\beta
    H_{\lfloor\beta^{-\frac 23}T\rfloor}\Bigr)\stackrel{\beta\da 0}
    {\Longrightarrow}\bigl((B^\sigma_t)_{t\in[0,T]},\widehat H^\sigma_T\bigr),
    \qquad \sigma,T>0.
    \end{equation}
As was shown in \cite{CR83}, the discrete local times process converges weakly
to the continuous local times process:
    \begin{equation}
    \label{betatozeroWC}
    \Bigl(\beta^{\frac 13}\ell_{\lfloor\beta^{-\frac 23}T\rfloor}
    (\lfloor x \beta^{-\frac 13}\rfloor)\Bigr)_{x\in \R}\stackrel{\beta\da 0}
    {\Longrightarrow}\big(L^\sigma(T,x)\big)_{x\in \R},\qquad \sigma,T>0.
    \end{equation}
This explains the scaling of the second component in \eqref{ntoinfty}--\eqref{betatozero}.
Since $(B^\sigma_t)_{t\geq 0}\stackrel{\Dcal}{=}(\sigma B_t)_{t\geq 0}$, we have that
        \begin{equation}
        \label{Bscaling}
        \big(L^\sigma(T,x)\big)_{x\in \R}\stackrel{\Dcal}{=}
        \big(\textstyle {\frac 1\sigma}
         L\big(T,\textstyle {\frac x\sigma})\big)_{x\in \R},
        \qquad \widehat H^\sigma_T\stackrel{\Dcal}{=}
        \textstyle {\frac 1\sigma}\widehat H_T, \qquad \sigma,T>0.
        \end{equation}

\medskip\noindent
II.\ The second invariance principle we will rely on was shown in \cite[Theorem 1.8]{Al},
and states that
    \begin{equation}
    \label{rtoinfty}
    \Bigl(\sigma^{-\frac 43}\bigl(\Ssig_{\lfloor\sigma^{\frac 23}t\rfloor}\bigr)_{t\in[0,T]},
    \1_{\bigl\{H_{\lfloor\sigma^{\frac 23}T\rfloor}=0\bigr\}}\Bigr)
    \stackrel{\sigma\to\infty}
    {\Longrightarrow}\bigl((B_t)_{t\in[0,T]},\1_{\{U>T\}}\bigr),\qquad T>0,
    \end{equation}
where the law of the random variable $U$ is given by its conditional distribution
given the underlying Brownian motion as
    \begin{equation}
    \widehat P\big(U>T\big|(B_t)_{t\in[0,T]}\big) = e^{-\widehat H_T},
    \end{equation}
and the limit $\sigma\to\infty$ is to be taken subject to the
following three technical restrictions:
    \begin{equation}
    \label{sigmares}
    \begin{array}{ll}
    &\mbox{(a) } \lim\limits_{N\to\infty}\limsup\limits_{\sigma\to\infty}
    E\bigl((\Ssig_1/\sigma)^2 \1_{\{|\Ssig_1/\sigma|>N\}}\bigr)=0;\\[0.3cm]
    &\mbox{(b) } \lim\limits_{\sigma\to\infty} \sigma^{\frac 23}
    \max\limits_{x\in\Z}P(\Ssig_1=x)=0;\\[0.3cm]
    &\mbox{(c) } \min\limits_{\sigma \geq 1} \min\limits_{0<|x|\leq c_1\sigma}
    \sigma P(\Ssig_1=x) \geq c_2 \mbox{ for some } c_1,c_2>0.
    \end{array}
    \end{equation}

The analogue of (\ref{betatozeroWC}) for $\sigma\to\infty$ under (\ref{sigmares})
is not known. Therefore, on top of (\ref{sigmares}), we will require a uniform
exponential moment for $\Ssig_1/\sigma$, i.e.,
    \begin{equation}
    \label{expmomsigma}
    \sup_{\sigma\geq 1}E(e^{\eps |\Ssig_1|/\sigma})<\infty
    \quad \mbox{for some }\eps>0,
    \end{equation}
which is obviously stronger than (\ref{sigmares})(a) and replaces the second condition
in \eqref{mv}. Note that the random walk with $P(S_1=x)=\frac{1}{2L}$ for
$x\in\{-L,\dots,-1,1,\dots,L\}$ satisfies \eqref{sigmares}--\eqref{expmomsigma}
(for which $\sigma^2\sim L^2/3$). So does the random walk with $P(S_1=x)=\frac{1}{2L}
(\frac{L-1}{L})^{|x|-1}$ for $x\in\Z\setminus\{0\}$ (for which $\sigma^2\sim L^2$).


\section{Large deviations}
\label{sec-LD}

To prove the results in Sections~\ref{mainres}--\ref{twodim}, we will actually
prove something much stronger, namely, {\it scaling of the large deviation rate
function for the empirical drift of the path}. We will show that the rate function
for the Domb-Joyce model and its variants scales to the rate function for the
Edwards model. Now, the existence of the rate function for the Domb-Joyce model
has been established only in a rather limited number of cases, namely, under
the assumption in \eqref{rcond}. In Section~\ref{sec-LDDJ} we summarize what
is known for this special case. For the variants of the Domb-Joyce model the
existence is still open. Therefore we will have to work with liminf's and limsup's.
The existence of the rate function for the Edwards model has been proved in our
recent paper \cite{HHK3} and its properties will be described in Section~\ref{sec-LDEM}.
Another important object is the cumulant generating function for the Edwards model,
which will be introduced in Section~\ref{sec-cum}. More refined large deviation
properties for the Edwards model also proved in \cite{HHK3}, which will be needed
in our proofs, are presented in Section~\ref{sec-refined}.

\subsection{Large deviations for the Domb-Joyce model.}
\label{sec-LDDJ}

Throughout this section we assume \eqref{rcond}. The main object of interest
in this section is the rate function $I_\beta$ defined by
\footnote{In fact, $I_\beta$ differs by a constant from what is usually called
a rate function: $I_\beta-r^*$ is the true rate function (see Theorem
\ref{thm-DJ*}(ii)).}
    \begin{equation}
    \label{rate}
    I_{\beta}(\theta) = -\lim_{n\to\infty}\frac
    1n\log E\bigl( e^{-\beta H_n} \1_{\{S_n\approx\theta n\}}\bigr)
    = -\lim_{n\to\infty}\frac
    1n\log\bigl\{Z_n^\beta Q_n^\beta(S_n\approx\theta n)\bigr\},
    \qquad \theta\in \R,
    \end{equation}
where $S_n\approx\theta n$ means that either $S_n=\lfloor\theta n\rfloor$
or $S_n=\lceil\theta n\rceil$ (possibly depending on the parity of these numbers).
For $\beta=\infty$ we adopt the convention $e^{-\infty H_n}=\1_{\{H_n=0\}}$.
Obviously, $I_{\beta}(\theta)=I_{\beta}(-\theta)$, and $I_{\beta}(\theta)
=\infty$ when $\theta>L$. Therefore we may restrict ourselves to $\theta\in [0,L]$.

Recall the three quantities $r^*,\theta^*,\sigma^*$ in Theorem~\ref{thm-DJ*}.
In the next theorem a fourth quantity $\theta^{**}$ appears, which, like the others,
depends on $\beta$ and on the distribution of $S_1$.

    \medskip
    \begin{theorem}[LDP]
    \label{thm-DJ}
    Fix $\beta\in [0,\infty]$, assume~\eqref{rcond}, and exclude the trivial
    case $(\beta,L)=(\infty,1)$.
    \begin{enumerate}
    \item[(i)]
    For any $\theta\in[0,L]$, the limit $I_{\beta}(\theta)$ in \eqref{rate} exists
    and is finite.
    \item[(ii)]
    $I_{\beta}$ is continuous and convex on $[0,L]$, and continuously
    differentiable on $(0,L)$.
    \item[(iii)]
    There is a number $\theta^{**}\in(0,\theta^*)$ such that $I_{\beta}$ is
    linearly decreasing on $[0,\theta^{**}]$, real-analytic and strictly convex
    on $(\theta^{**},L)$, and attains its unique minimum at $\theta^*$ with
    height $I_{\beta}(\theta^*)=\r^*$ and curvature $I''_{\beta}(\theta^*)
    =1/\sigma^{*2}$.
    \end{enumerate}
    \end{theorem}


\vskip 1truecm

\setlength{\unitlength}{0.7cm}

\begin{picture}(12,8)(-4,-1.5)

  \put(0,-2){\line(12,0){12}}
  \put(0,-2){\line(0,8){8}}
  \put(-.1,-3){$0$}
  {\thicklines
   \qbezier(0,3)(2,2.5)(4,2)
   \qbezier(4,2)(6,1.55)(7,1.85)
   \qbezier(7,1.85)(8,2.2)(9,5)
  }
  \qbezier[40](4,-2)(4,0)(4,2)
  \qbezier[40](6,-2)(6,-.25)(6,1.75)
  \qbezier[70](0,1.7)(4,1.7)(8,1.7)
  \qbezier[70](4,2)(6,1.5)(7,1.25)
  \qbezier[70](9,-2)(9,1.5)(9,5)
  \put(4,2){\circle*{.15}}
  \put(0,3){\circle*{.15}}
  \put(6,1.74){\circle*{.15}}
  \put(9,5){\circle*{.15}}
  \put(3.4,-3){$\theta^{**}(\beta)$}
  \put(5.5,-3){$\theta^*(\beta)$}
  \put(8.9,-3){$L$}
  \put(12.5,-2.15){$\theta$}
  \put(-1.7,1.6){$\r^*(\beta)$}
  \put(-.3,6.5){$I_{\beta}(\theta)$}
  \put(1,-4.3){\small
  Fig.\ 1. Qualitative picture of $\theta\mapsto I_{\beta}(\theta)$.
  \normalsize}
  \end{picture}

\vskip 2.5truecm


\noindent
Theorem~\ref{thm-DJ} is proved for simple random walk ($L=1$) in \cite[Theorem IX.32]{dHLD},
relying on the methods and results of \cite{GH}. We have checked that this proof can
be extended to general $L\in\N$ with the help of the methods and results of \cite{K2}.

The main ingredients of the proof of Theorem~\ref{thm-DJ} are reflection arguments
and precise analytic knowledge of the contribution to the intersection local
time coming from paths that satisfy the so-called ``bridge condition'', i.e., lie
between their starting and ending locations $S_0$ and $S_n$. The linear piece of
the rate function has the following intuitive explanation. If $\theta\geq\theta^{**}$,
then the optimal strategy for the path to realize $S_n\approx\theta n$ is to assume local
drift $\theta$ during $n$ steps. In particular, the path then satisfies the bridge
condition, and this reasoning leads to the strict convexity and real-analyticity of
the rate function on $(\theta^{**},L)$. If, on the other hand, $0\leq\theta<\theta^{**}$,
then this strategy is too expensive, since too small a drift leads to too many
self-intersections. Therefore the optimal strategy now is to move with local drift
$\theta^{**}$ during $\frac{\theta^{**}+\theta}{2\theta^{**}}n$ steps and with
local drift $-\theta^{**}$ during the remaining $\frac{\theta^{**}-\theta}{2\theta^{**}}n$
steps, thus making an overshoot of size $\frac{\theta^{**}-\theta}{2}n$, and this
reasoning leads to the linearity of the rate function on $[0,\theta^{**}]$.

\subsection{Large deviations for the Edwards model.}
\label{sec-LDEM}

The analogue of \eqref{rate} for the Edwards model is the rate function
$\widehat I_\beta$ defined by
    \begin{equation}
    \label{rateED}
    \widehat I_\beta(b) = - \lim_{T\to\infty}\frac 1T\log
    \widehat E\Bigl(e^{-\beta\widehat H_T} \1_{\{B_T\approx bT\}}\Bigr)
    = -\lim_{T\to\infty}\frac 1T\log
    \left\{\widehat Z^\beta_T \widehat Q_T^\beta\left(B_T\approx bT\right)\right\},
    \qquad b\in\R,
    \end{equation}
where $B_T \approx bT$ means that $|B_T-bT|\leq \gamma_T$ for some $\gamma_T>0$ such
that $\gamma_T/T\to0$ and $\gamma_T/\sqrt T\to\infty$ as $T\to\infty$. In \cite{HHK3}
we proved that the limit in \eqref{rateED} exists and is independent of the choice
of $\gamma_T$. From (\ref{scalEM}) it is clear that this rate function satisfies
the scaling relation
\begin{equation}
\label{betascalrat}
\beta^{-\frac 23}\widehat I_\beta(\beta^{\frac 13}\cdot)
=\widehat I_1(\cdot),
\end{equation}
provided the limit in \eqref{rateED} exists for $\beta=1$.

Recall the three quantities $a^*,b^*,c^*$ in Theorem~\ref{thm-CLTEDW}. In the next
theorem a fourth quantity $b^{**}$ appears.

    \medskip
    \begin{theorem}[LDP]
    \label{thm-LDEM}
    $\mbox{}$
    \begin{enumerate}
    \item[(i)]
    For any $b\in[0,\infty)$, the limit $\widehat I_1(b)$ in \eqref{rateED}
    exists and is finite (and is independent of the choice of $\gamma_T$).
    \item[(ii)]
    $\widehat I_1$ is continuous and convex on $[0,\infty)$, and continuously
    differentiable on $(0,\infty)$.
    \item[(iii)]
    There is a number $b^{**}\in(0,b^*)$ such that $\widehat I_1$
    is linearly decreasing on $[0,b^{**}]$, real-analytic and strictly convex on
    $(b^{**},\infty)$, and attains its unique minimum at $b^*$ with
    height $\widehat I_1(b^*)=a^*$ and curvature $\widehat I_1''(b^*)
    =1/c^{*2}$.
    \end{enumerate}
    \end{theorem}
\medskip
\noindent
Theorem~\ref{thm-LDEM} is proved in \cite{HHK3}. The numerical value of $b^{**}$
is $b^{**}\approx 0.85$. Note the close analogy with Theorem~\ref{thm-DJ}. The
linear piece has the same intuitive explanation in terms of overshoot.


\vskip 1truecm

\setlength{\unitlength}{0.7cm}

\begin{picture}(12,8)(-4,-1.5)

  \put(0,-2){\line(12,0){12}}
  \put(0,-2){\line(0,8){8}}
  \put(-.1,-3){$0$}
  {\thicklines
   \qbezier(0,3)(2,2.5)(4,2)
   \qbezier(4,2)(6,1.55)(7,1.85)
   \qbezier(7,1.85)(8,2.2)(9,5)
  }
  \qbezier[40](4,-2)(4,0)(4,2)
  \qbezier[40](6,-2)(6,-.25)(6,1.75)
  \qbezier[70](0,1.7)(4,1.7)(8,1.7)
  \qbezier[70](4,2)(6,1.5)(7,1.25)
  \qbezier[10](9,5)(9.15,5.4)(9.3,6.2)

  \put(4,2){\circle*{.15}}
  \put(6,1.74){\circle*{.15}}
  \put(0,3){\circle*{.15}}
  \put(3.8,-3){$b^{**}$}
  \put(5.9,-3){$b^*$}
  \put(12.5,-2.15){$b$}
  \put(-1,1.62){$a^*$}
  \put(-.25,6.5){$\widehat I_1(b)$}
  \put(1,-4.2){\small
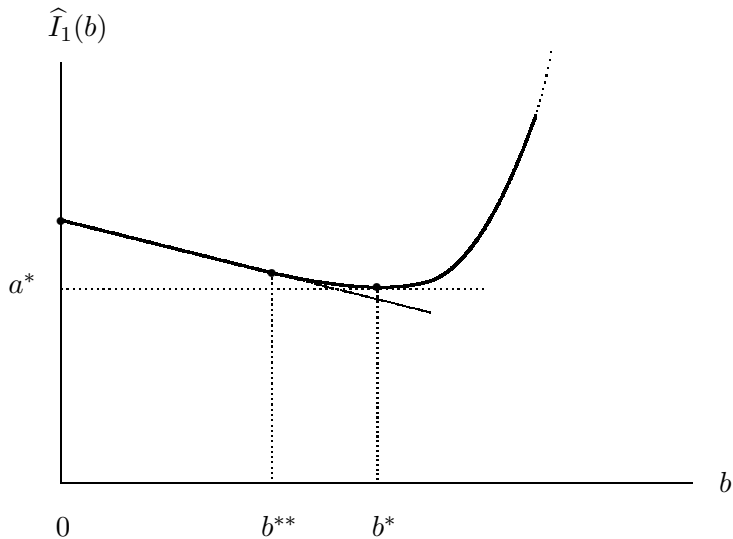
  Fig.\ 2. Qualitative picture of $b\mapsto \widehat I_1(b)$.
  \normalsize}
  \end{picture}

\vskip 2.5truecm


Denote by $\widehat I_\beta^\sigma$ the rate function in \eqref{rateED} for
the Brownian motion with generator $\frac 12\sigma^2\Delta$. Like $\widehat I_\beta$,
it satisfies the scaling relation $\beta^{-\frac 23}\widehat I_\beta^\sigma
(\beta^{\frac 13}\cdot)=\widehat I^\sigma_1(\cdot)$ in \eqref{betascalrat}.
Furthermore, from \eqref{Bscaling} we obtain the scaling relation
        \begin{equation}
        \label{sigmascal}
        \widehat I^\sigma_1(\cdot)
        = \sigma^{-\frac 23}\widehat I_1(\sigma^{-\frac 23}\cdot).
        \end{equation}

\subsection{Cumulant generating function for the Edwards model.}
\label{sec-cum}

There is an intimate connection between the rate function in \eqref{rateED}
and the cumulant generating function $\Lambda^+\colon \R\to\R$ given by
   \begin{equation}
   \label{Lambdadef}
   \Lambda^+(\mu)=\lim_{T\to\infty}\frac 1T\log
   \widehat E\bigl(e^{-\widehat H_T} e^{\mu B_T}\1_{\{B_T\geq 0\}}\bigr),\qquad
   \mu\in\R.
\end{equation}

   \begin{prop}[Exponential moments]
   \label{exponmom}
   $\mbox{}$
   \begin{enumerate}
   \item[(i)]
   For any $\mu\in\R$, the limit $\Lambda^+(\mu)$ in \eqref{Lambdadef} exists
   and is finite.
   \item[(ii)]
   There is a number $\rho(a^{**})>0$ such that $\Lambda^+$ is constant on
   $(-\infty,-\rho(a^{**})]$, and strictly increasing, strictly convex and
   real-analytic on $(-\rho(a^{**}),\infty)$. In $-\rho(a^{**})$, $\Lambda^+$
   is continuous, but not differentiable.
   \item[(iii)]
   $\lim_{\mu\da-\rho(a^{**})}(\Lambda^+)'(\mu)=b^{**}$, $(\Lambda^+)'(0)=b^*$,
   and $\lim_{\mu\to\infty}(\Lambda^+)'(\mu)=\infty$.
   \item[(iv)]
   The restriction of $\widehat I_1$ to $[0,\infty)$ is the Legendre transform of
   $\Lambda^+$, i.e.,
   \begin{equation}
   \label{LegTra}
   \widehat I_1(b)= \max_{\mu\in \R}\bigl[\mu b-\Lambda^+(\mu)\bigr],
   \qquad b\geq 0.
   \end{equation}
   \end{enumerate}
   \end{prop}

\medskip\noindent
Proposition~\ref{exponmom} is proved in \cite{HHK3}. The numerical value of
$\rho(a^{**})$ is $\rho(a^{**})\approx 0.78$. By \eqref{LegTra}, $-\rho(a^{**})$
is the slope of the linear piece in Fig.\ 2. Note that $\Lambda^+(0)=-a^*$ by
Theorem~\ref{thm-CLTEDW}(ii) and \eqref{betascalrat}.

As a consequence of Proposition~\ref{exponmom}(ii), the maximum on the right-hand
side of \eqref{LegTra} is attained in some $\mu>-\rho(a^{**})$ if $b>b^{**}$ and
in $\mu=-\rho(a^{**})$ if $0\leq b\leq b^{**}$.

Let $\Lambda^-$ denote the cumulant generating function with $\1_{\{B_T\leq 0\}}$
instead of $\1_{\{B_T\geq 0\}}$. Then analogous assertions for $\Lambda^-$
hold as well. In particular, the restriction of $\widehat I_1$ to $(-\infty,0]$
is the Legendre transform of $\Lambda^-$. By symmetry, $\Lambda^+(-\mu) =\Lambda^-(\mu)$
for any $\mu\in\R$. Consequently, the cumulant generating function $\Lambda(\mu)
=\lim_{T\to\infty}\frac 1T\log \widehat E\bigl(e^{-\widehat H_T}e^{\mu B_T}\bigr)
=\Lambda^+(\mu)\vee\Lambda^-(\mu)=\Lambda^+(|\mu|)$ exists for any $\mu\in\R$
and is not differentiable at 0.

Let $\Lambda^+_{\sigma}$ and $\Lambda^-_{\sigma}$ denote the corresponding cumulant
generating functions for the Edwards model with variance $\sigma^2$ (i.e., where
the generator of the underlying Brownian motion is $\frac 12\sigma^2\Delta$). Then
we have the scaling relation $\sigma^{\frac 23}\Lambda_{\sigma}^+(\sigma^{-\frac 43}\,
\cdot)=\Lambda^+(\cdot)$. Moreover, we have
        \begin{equation}
        \label{LegTra2}
        \widehat I_1^{\sigma}(b)
        = \max_{\mu\in\R}\bigl[\mu b-\Lambda_{\sigma}^+(\mu)\bigr]
        =\begin{cases}\max\limits_{\mu\geq 0}
        \bigl[\mu b-\Lambda_{\sigma}^+(\mu)\bigr]&\mbox{if }b\geq b^*\sigma^{\frac 23},\\
        \max\limits_{\mu\leq 0}\bigl[\mu b-\Lambda_{\sigma}^+(\mu)\bigr]&\mbox{if }
        0\leq b\leq b^*\sigma^{\frac 23}.
        \end{cases}
        \end{equation}
Analogous assertions hold for $\Lambda^-_{\sigma}$.

\subsection{More refined large deviation properties for the Edwards model.}
\label{sec-refined}

In the proofs we will need some further refinements of Proposition~\ref{exponmom}.
Abbreviate $B_{[0,T]}=(B_t)_{t\in[0,T]}$. For $T>0$, $\delta,C\in(0,\infty]$
and $\alpha\in[0,\infty)$, define events
    \begin{eqnarray}
    \widehat\ev(\delta,T)
    &=& \bigl\{B_{[0,T]}\subset[-\delta,B_T+\delta]\bigr\},
    \label{Ehatdef}\\
    \widehat \ev^{\leq}(\delta,C;T)
    &=&\Big\{\max_{x\in[-\delta,\delta]}
    L(T,x)\leq C,\max_{x\in[B_T-\delta,B_T+\delta]}L(T,x)\leq C\Big\},\label{Fhatdeflq}\\
    \widehat \ev^{\,\geq}(\delta,\alpha;T)
    &=&\Big\{\max_{x\in[B_T-\delta,B_T+\delta]}
    L(T,x)\geq\alpha\delta^{-\frac 12}\Big\}.\label{Fhatdefgq}
    \end{eqnarray}
Note that $\widehat\ev^{\leq}(\delta,\infty;T)$ and $\widehat\ev^{\,\geq}(\delta,0;T)$
are the full space.

\newpage

    \medskip\noindent
    \begin{prop}[Overshoots]
    \label{prop-ingred}
    Fix $\mu>-\rho(a^{**})$. Then:
    \begin{enumerate}
    \item[(i)]
    For any $\delta,C\in(0,\infty]$ there exists a $K_1(\delta,C)\in(0,\infty)$
    such that
    \begin{equation}
    \label{refine1}
    e^{-\Lambda^+(\mu)T}\widehat E\Bigl(e^{-\widehat H_T}e^{\mu B_T}\1_{\widehat \ev(\delta,T)}
    \1_{\widehat \ev^\leq(\delta,C;T)}\1_{\{B_T\geq 0\}}\Bigr) = K_1(\delta,C)+o(1),
    \qquad T\to\infty.
    \end{equation}
    Moreover, if $\mu=\mu_b$ solves $\widehat I_1(b)=\mu b-\Lambda^+(\mu)$, then
    the same is true when $\1_{\{B_T\geq 0\}}$ is replaced by $\1_{\{B_T\approx bT\}}$.
    \item[(ii)]
    For any $\delta,\alpha\in(0,\infty)$ there exists a $K_2(\delta,\alpha)\in(0,\infty)$
    such that
    \begin{equation}
    \label{refine2}
    e^{-\Lambda^+(\mu)T}\widehat E\Bigl(e^{-\widehat H_T}e^{\mu B_T}\1_{\widehat \ev(\delta,T)}
    \1_{\widehat \ev^{\,\geq}(\delta,\alpha;T)}\1_{\{B_T\geq 0\}}\Bigr)
    = K_2(\delta,\alpha)+o(1), \qquad T\to\infty.
    \end{equation}
    \item[(iii)] For any $\alpha\in(0,\infty)$,
    \begin{equation}\label{Kbound}
    \lim_{\delta\downarrow 0}\frac{K_2(\delta,\alpha)}{K_1(\delta,\infty)}=0.
    \end{equation}
    \end{enumerate}
    \end{prop}

\medskip\noindent
Proposition~\ref{prop-ingred} is proved in \cite{HHK3}.


\section{Proof of Theorem \ref{ratescal}(i)}
\label{sec-proof(i)}

In this section we consider the limit $\beta\da 0$. Let $(S_n)_{n\in\N_0}$
be a random walk satisfying \eqref{mv}. As announced at the beginning of
Section~\ref{sec-LD}, we will identify the scaling limit of the entire
large deviation rate function (for the linear asymptotics of the endpoint)
for the Domb-Joyce model in terms of that for the Edwards model, and we will
deduce Theorem~\ref{ratescal}(i) from this scaling limit. However, as pointed
out at the beginning of Section~\ref{sec-LD}, the existence of the rate function
has not been established in full generality for the Domb-Joyce model, and we
will make no attempt to do so. Instead, we will be working with {\it approximative
rate functions}, which are defined as a limsup or a liminf instead of a lim.

\subsection{Approximative large deviations.}
\label{sec-appro}

It will be sufficient to deal with the event $\{S_n\geq \theta n\}$ for $\theta$
to the right of the scaled minimum point of the limiting rate function, and with
$\{S_n\leq \theta n\}$ for $\theta$ to the left of it. To this end, define
    \begin{equation}
    \label{rate+}
    I^+_{\beta}(\theta;\widetilde\theta)=
    \begin{cases}
    -\liminf\limits_{n\to\infty}\frac
    1n\log E\bigl( e^{-\beta H_n} \1_{\{\Ssig_n\geq \theta n\}}\bigr)
    &\text{ if }\theta\geq\widetilde\theta,\\
    -\liminf\limits_{n\to\infty}\frac
    1n\log E\bigl( e^{-\beta H_n} \1_{\{0\leq\Ssig_n\leq \theta n\}}\bigr)
    &\text{ if }\theta\leq\widetilde \theta,
    \end{cases}
    \end{equation}
and define $I^-_{\beta}(\theta;\widetilde\theta)$ in the same way with $\limsup$ instead
of $\liminf$. For $\beta=\infty$, recall the convention $e^{-\infty H_n}=\1_{\{H_n=0\}}$.

In the special case of \eqref{rcond}, we know from Theorem \ref{thm-DJ} that the
limit $I_\beta(\theta)$ in \eqref{rate} exists. Since $I_{\beta}$ is unimodal with
unique minimiser $\theta^*$, it follows that both limits in \eqref{rate+} exist
and that
    \begin{equation}
    I^+_{\beta}(\theta;\theta^*)=I^-_{\beta}(\theta;\theta^*)
    =I_{\beta}(\theta),\qquad 0 \leq\theta\leq L.
    \end{equation}

Our main result in this section shows that the approximative rate function in
\eqref{rate+} scales, as $\beta \downarrow 0$, to the rate function for the
Edwards model with parameter $\sigma$.

    \medskip
    \begin{prop}
    \label{ratescal2}
    Fix $\sigma\in(0,\infty)$. Then, under \eqref{mv},
     \begin{eqnarray}
     \liminf_{\beta\da0}\beta^{-\frac 23}
     I^-_{\beta}\bigl(b \beta^{\frac 13};b^*\beta^{\frac 13}\sigma^{\frac 23}\bigr)
     &\geq&\widehat I^\sigma_1(b),\qquad b\geq 0,\label{betascal}\\
     \limsup_{\beta\da0}\beta^{-\frac 23}I^+_{\beta}
     \bigl(b\beta^{\frac 13};b^*\beta^{\frac 13}\sigma^{\frac 23}\bigr)
     &\leq&\widehat I^\sigma_1(b),\qquad b>b^{**}\sigma^{\frac 23}.
     \label{betascallow}
     \end{eqnarray}
     \end{prop}

\medskip\noindent
Proposition~\ref{ratescal2} is proved in Section~\ref{sec-ratescal2}. In the
special case of (\ref{rcond}), we infer from Theorem~\ref{thm-DJ} and
Proposition~\ref{ratescal2} that
    \begin{equation}
    \label{ratescalperf}
    \lim_{\beta\da0}\beta^{-\frac 23}
    I_{\beta}\bigl(b \beta^{\frac 13}\bigr)
    =\widehat I^\sigma_1(b),\qquad b>b^{**}\sigma^{\frac 23}.
    \end{equation}

\medskip\noindent
The proof of \eqref{betascallow} for $0\leq b\leq b^{**}\sigma^{\frac 23}$
remains open. To extend \eqref{betascallow} to this regime would require
some further refinements of our method (see Section~\ref{sec-disc}).

\subsection{Proof of Theorem \ref{ratescal}(i) and Corollary~\ref{asy}.}
\label{sec-proofcor1}

\medskip\noindent
{\bf 1.} Fix $\eps>0$. We will show that, for $\beta>0$ sufficiently small,
    \begin{equation}
    \label{LLNupper}
    \lim_{n\to\infty}\frac 1n\log Q_n^\beta\Bigl(\frac{|\Ssig_n|}
    {\beta^{\frac 13}n}-b^*\sigma^{\frac 23}>\eps\Bigr)<0.
    \end{equation}
This obviously implies the upper half of the statement in \eqref{LLN}.
The lower half can be derived in the same manner.

\medskip\noindent
{\bf 2.} To prove \eqref{LLNupper}, put $b'=b^*\sigma^{\frac 23}+\frac\eps2$
and $b=b^*\sigma^{\frac 23}+\eps$. Since $\widehat I^\sigma_1$ is strictly
increasing on $[b^*\sigma^{\frac 23},\infty)$, it is possible to pick $\gamma>0$
so small (depending on $\eps$) that
    \begin{equation}
    \label{gammapick}
    \widehat I^\sigma_1(b)-\widehat I^\sigma_1(b')-2\gamma>0.
    \end{equation}
According to Proposition~\ref{ratescal2}, we may pick $\beta>0$ so small
(depending on $\gamma$) that
    \begin{equation}
    \label{betapick}
    I^-_{\beta}\bigl(b \beta^{\frac 13};b^*\beta^{\frac 13}\sigma^{\frac 23}\bigr)
    \geq\bigl[\widehat I^\sigma_1(b)-\gamma\bigr]\beta^{\frac 23},
    \qquad
    I^+_{\beta}\bigl(b'\beta^{\frac 13};b^*\beta^{\frac 13}\sigma^{\frac 23}\bigr)
    \leq\bigl[\widehat I^\sigma_1(b')+\gamma\bigr]\beta^{\frac 23}.
    \end{equation}
Now we can bound (recall \eqref{PM})
    \begin{equation}
    \label{LLNesti}
    \begin{aligned}
    Q_n^\beta\Bigl(\frac{\Ssig_n}{\beta^{\frac 13}n}-b^*\sigma^{\frac 23}>\eps\Bigr)
    &=\frac{E\bigl(e^{-\beta H_n}\1_{\{\Ssig_n>b\beta^{\frac 13}n\}}\bigr)}
    {E\bigl(e^{-\beta H_n}\bigr)}
    \leq \frac{E\bigl(e^{-\beta H_n}\1_{\{\Ssig_n>b\beta^{\frac 13}n\}}\bigr)}
    {E\bigl(e^{-\beta H_n}\1_{\{\Ssig_n>b'\beta^{\frac 13}n\}}\bigr)}\\
    &\leq \exp\Bigl\{-n\bigl[I^-_{\beta}\bigl(b \beta^{\frac 13}
    ;b^*\beta^{\frac 13}\sigma^{\frac 23}\bigr)-
    I^+_{\beta}\bigl(b'\beta^{\frac 13};b^*\beta^{\frac 13}\sigma^{\frac 23}\bigr)\bigr]
    +o(n)\Bigr\},
    \end{aligned}
    \end{equation}
where we use the definitions of $I^-_{\beta}$ and $I^+_{\beta}$. Insert
\eqref{gammapick}--\eqref{betapick}, to see that the term between square
brackets in the exponent of \eqref{LLNesti} is strictly positive. This implies
\eqref{LLNupper}.

\medskip\noindent
{\bf 3.} The proof of Corollary~\ref{asy} is as follows. Assume \eqref{rcond}.
First, by \eqref{ratescalperf}, the function $f_\beta$ defined by $f_\beta(\cdot)=
\beta^{-\frac 23} I_{\beta}\bigl(\beta^{\frac 13}\cdot\bigr)$ converges to
$\widehat I_1^\sigma$ on $(b^{**}\sigma^{\frac 23},\infty)$. In particular,
the unique minimal value of $f_\beta$, which is $r^*(\beta)\beta^{-\frac 23}$ by
Theorem~\ref{thm-DJ}, converges to the unique minimal value of $\widehat I_1^\sigma$,
which is $a^*\sigma^{-\frac 23}$ by Theorem~\ref{thm-LDEM}. This proves the
first assertion in \eqref{asyvalue}. Next, by \eqref{ratescalperf}, $f_\beta$
converges to $\widehat I_1^{\sigma}$ in the three points $b^*\sigma^{\frac 23}-\eps$,
$b^*\sigma^{\frac 23}$ and $b^*\sigma^{\frac 23}+\eps$ for $\eps>0$ small enough.
For $\beta$ small enough, both $f_\beta(b^*\sigma^{\frac 23}-\eps)$ and
$f_\beta(b^*\sigma^{\frac 23}+\eps)$ are strictly larger than $f_\beta(b^*\sigma^{\frac 23})$.
By unimodality, this implies that the unique minimiser of $f_\beta$, which is
$\theta^*(\beta)\beta^{-\frac 13}$ by Theorem~\ref{thm-DJ}, lies in
$(b^*\sigma^{\frac 23}-\eps,b^*\sigma^{\frac 23}+\eps)$. Let $\eps\da 0$ to
obtain the second assertion in \eqref{asyvalue}.
\qed

Note that convexity of $f_\beta$ yields that even $(f_\beta)'$ converges to
$(\widehat I_1^\sigma)'$. However, we have no control over $(f_\beta)''$,
which is why we are unable to prove \eqref{sigscal}.

\subsection{Proof of Proposition~\ref{ratescal2}.}
\label{sec-ratescal2}

In Section~\ref{upperbound} we prove \eqref{betascal}, in Section~\ref{lowerbound*}
we prove \eqref{betascallow}. The main idea is to cut the path into smaller pieces
to which the weak convergence assertion in \eqref{betatozero} can be applied. The
mutual interaction between the pieces has to be controlled appropriately. This
is done by providing estimates in which either the pieces are independent or there
is an interaction only between neighboring pieces. We define
    \begin{equation}
    H_n'=\sum_{\heap{i,j=1}{i\not= j}}^n
    \1_{\{S_i=S_j\}}=H_n -2(\ell_n(0)-1).
    \label{Hn'def}
    \end{equation}
The proof runs via the moment
generating function
    \begin{equation}
    \label{momZ}
    Z_n^{\beta}(\mu) = E\bigl(e^{-\beta H_n'} e^{\mu \beta^{\frac 13}
    S_n}\bigr),\qquad n\in\N, \mu\in\R
    \end{equation}
which is the discrete analogue of the expectation in \eqref{Lambdadef}.

\subsubsection{Proof of (\ref{betascal})}
\label{upperbound}

\noindent
{\bf 1.}\ Fix $b\geq b^*\sigma^{\frac 23}$. Use the exponential Chebyshev inequality
to get the following upper bound for $\mu\geq 0$:
    \begin{equation}
    \label{esti1}
    E\bigl(e^{-\beta H_n}\1_{\{\Ssig_n \geq b \beta^{\frac 13}n
    \}}\bigr) \leq e^{-\mu b \beta^{\frac 23}n} Z_n^{\beta}(\mu).
    \end{equation}
Fix a large auxiliary parameter $T>0$ and abbreviate $T_\beta =
\beta^{-\frac 23}T$. Split the path of length $n$ into $n/T_\beta$
pieces of length $T_\beta$. (To simplify the notation, assume that
both $n/T_\beta$ and $T_\beta$ are integers.) Drop the interaction
between any two of the pieces, to obtain an upper bound on
$Z_n^\beta(\mu)$. After the pieces are decoupled they are independent
of each other. This reasoning yields
    \begin{equation}
    \label{esti1a}
    Z_n^{\beta}(\mu)
    \leq \bigl(Z_{T_\beta}^{\beta}(\mu)\bigr)^{n/T_\beta}.
    \end{equation}
Substitute this estimate into \eqref{esti1}, take logs, divide by
$\beta^{{\frac 23}}n$ and let $n\to\infty$, to obtain (recall \eqref{rate+})
    \begin{equation}
    \label{esti2}
     \begin{aligned}
    \beta^{-\frac 23}I_{\beta}^-(b\beta^{\frac 13}
    ;b^*\beta^{\frac 13}\sigma^{\frac 23})
    &\geq -\beta^{-\frac 23}\liminf_{n\to\infty}
    \frac 1n\log\bigl(\mbox{l.h.s.~of \eqref{esti1}}\bigr)\\
    &\geq -\beta^{-\frac 23}\liminf_{n\to\infty}
    \frac 1n\log\Bigl[e^{-\mu b \beta^{\frac 23}n}
    \bigl(Z_{T_\beta}^{\beta}(\mu)\bigr)^{n\beta^{\frac 23}/T}\Bigr]\\
    &=\mu b -\frac 1T\log Z_{T_\beta}^{\beta}(\mu).
    \end{aligned}
    \end{equation}

\medskip\noindent
{\bf 2.}\ The next lemma states that, under \eqref{mv}, the expectation in the
right-hand side of \eqref{esti2} converges to the corresponding Brownian
expectation. Its proof is given in part 4.

    \medskip
    \begin{lemma}
    \label{lem-wcbeta}
    Assume \eqref{mv}. Then, for any $\mu\in\R$,
        \begin{equation}
        \lim_{\beta\da 0} Z_{T_\beta}^{\beta}(\mu)
        =\widehat E(e^{-\widehat H_T^{\sigma}} e^{\mu B_T^{\sigma}}).
        \end{equation}
    \end{lemma}

\medskip\noindent
Lemma~\ref{lem-wcbeta} applied to \eqref{esti2} yields
    \begin{equation}
    \liminf_{\beta\downarrow 0} [\beta^{-\frac 23}
    I_{\beta}^-(b\beta^{\frac 13};b^*\beta^{\frac 13}\sigma^{\frac 23})]
    \geq \mu b-\frac1T \log\widehat E(e^{-\widehat H_T^{\sigma}} e^{\mu B_T^{\sigma}}),
    \qquad \mu\geq 0.
    \end{equation}
Now let $T\to\infty$ and use \eqref{Lambdadef}, to obtain
    \begin{equation}
    \label{esti3}
    \liminf_{\beta\downarrow 0}
    [\beta^{-\frac 23}I_{\beta}^-(b\beta^{\frac 13}
    ;b^*\beta^{\frac 13}\sigma^{\frac 23})]
    \geq \mu b-\Lambda^+_{\sigma}(\mu).
    \end{equation}
Maximize over $\mu\geq 0$ and use (\ref{LegTra2}), to arrive at the assertion in
\eqref{betascal}.

\medskip\noindent
{\bf 3.}\ The proof for $0\leq b\leq b^*\sigma^{\frac 23}$ follows the same pattern.
Estimate, for $\mu\leq 0$,
    \begin{equation}
    E\Bigl(e^{-\beta H_n}\1_{\{0\leq\Ssig_n\leq b\beta^{\frac 13}n\}}\Bigr)
    \leq e^{-\mu b \beta^{\frac 23}n} Z_n^{\beta}(\mu).
    \end{equation}
In the same way as above we obtain
    \begin{equation}
    \label{esti4}
    \liminf_{\beta \downarrow 0}
    [\beta^{-\frac 23}I_{\beta}^-(b\beta^{\frac 13};b^*\beta^{\frac 13}\sigma^{\frac 23})]
    \geq \mu b-\Lambda_\sigma^+(\mu).
    \end{equation}
Now maximize over $\mu\leq 0$ and again use \eqref{LegTra2}.
\qed

\medskip\noindent
{\bf 4.}\ We finish by proving Lemma \ref{lem-wcbeta}.

\begin{proofsect}{Proof of Lemma \ref{lem-wcbeta}.}
Fix $\mu\in\R$. By the weak convergence assertion in \eqref{betatozero}, together
with dominated convergence, we have for every $K>0$,
    \begin{equation}
    \label{Kineq}
    \lim_{\beta\da 0} E\Bigl(e^{-\beta H_{T_\beta}'}
    e^{\mu\beta^{\frac 13}\Ssig_{T_\beta}}
    \1_{\{\beta^{\frac 13}|\Ssig_{T_\beta}|<K\}}\Bigr)
    =\widehat E\left(e^{-\widehat H_T^{\sigma}}
    e^{\mu B_T^{\sigma}}\1_{\{|B_{T}^{\sigma}|<K\}}\right).
    \end{equation}
The right-hand side of \eqref{Kineq} increases to $\widehat E(e^{-\widehat
H_T^{\sigma}} e^{\mu B_T^{\sigma}})$ as $K\rightarrow \infty$.
Therefore it suffices to show that
    \begin{equation}
    \label{K2}
    \lim_{K\rightarrow \infty}\limsup_{\beta\da 0}
    E\Bigl(e^{-\beta H_{T_\beta}'}e^{\mu\beta^{\frac 13}\Ssig_{T_\beta}}
    \1_{\{\beta^{\frac 13}|\Ssig_{T_\beta}| \geq K\}}\Bigr)=0.
    \end{equation}
To prove \eqref{K2}, use the Cauchy-Schwarz inequality:
    \begin{equation}
    E\Bigl(e^{-\beta H_{T_\beta}'}e^{\mu\beta^{\frac 13}\Ssig_{T_\beta}}
    \1_{\{\beta^{\frac 13}|\Ssig_{T_\beta}| \geq K\}}\Bigr)^2\leq
    P\bigl(\beta^{\frac 13}|\Ssig_{T_\beta}| \geq K\bigr)
    E\bigl(e^{2\mu\beta^{\frac 13}\Ssig_{T_\beta}}\bigr).
    \end{equation}
The first term converges to $\widehat P(|B_T|\geq K)$ as $\beta\da 0$, which vanishes as
$K\to\infty$. Therefore it suffices to show that
    \begin{equation}
    \label{K3}
    \limsup_{\beta\downarrow 0} E\bigl(e^{2\mu\beta^{\frac 13}\Ssig_{T_\beta}}\bigr)<\infty.
    \end{equation}
To prove \eqref{K3}, denote the moment generating function of $\Ssig_1$ by $\varphi(t)
=E(e^{t\Ssig_1})$, $t\in\R$. Then
    \begin{equation}
    \label{lapo1}
    E(e^{2\mu\beta^{\frac 13}\Ssig_{T_\beta}})=\varphi(2\mu\beta^{\frac13})^{T_{\beta}}.
    \end{equation}
By \eqref{mv}, the right-hand side is finite for $\beta$ small enough (depending
on $\mu$). Now note that, by \eqref{mv}--\eqref{sigma},
    \begin{equation}
    \label{lapo2}
    \varphi(t)=1+\frac 12 \sigma^2 t^2 +{\Ocal}(|t|^3),\qquad t\to 0.
    \end{equation}
Put $t=2\mu\beta^{\frac13}$ and combine \eqref{lapo1}--\eqref{lapo2}, to get
    \begin{equation}
    E(e^{2\mu\beta^{\frac 13}\Ssig_{T_\beta}})\leq
    e^{T_\beta[\frac{1}{2}\sigma^2 t^2 +{\Ocal}(|t|^3)]}
    =e^{2\mu^2\sigma^2 T[1+{\Ocal}(\beta^{\frac13})]}, \qquad \beta\da 0.
    \end{equation}
This proves \eqref{K3} and completes the proof of Lemma~\ref{lem-wcbeta}.
\end{proofsect}
\qed

\subsubsection{Proof of (\ref{betascallow})}
\label{lowerbound*}

We again cut the path into pieces as in Section~\ref{upperbound}, but this time we
keep control of the interaction between the pieces. Since we are looking for
a lower bound on an expectation, we may freely require additional properties
of the pieces in such a way that we can control their mutual interaction and
still perform the limit $\beta\da0$.

\medskip\noindent
{\bf 1.}\ Fix $b\geq b^*\sigma^{\frac 23}$. We require that in each piece the
path has speed $\geq b\beta^{\frac 13}$, does not go too far beyond its starting
and ending locations, and has local times in the overlapping areas
that are uniformly bounded by a constant. To formulate this precisely,
for $i=1,\dots,n/T_\beta$ denote by
    \begin{equation}
    \label{piece}
    S^{\smallsup{i}} = (S_j^{\smallsup{i}})_{j=0}^{T_\beta}
    \quad \mbox{ with } \quad S_j^{\smallsup{i}}=S_{j+(i-1)T_\beta}-S_{(i-1)T_\beta}
    \end{equation}
the $i$-th piece shifted such that it starts at the origin, and denote by
    \begin{equation}
    \ell^{\smallsup{i}}(x) =\sum_{j=(i-1)T_\beta+1}^{iT_\beta}
    \1_{\{S_j-S_{(i-1)T_\beta}=x\}}
    =\sum_{j=1}^{T_\beta}\1_{\{S^{\smallsup{i}}_j=x\}},\qquad x\in\Z,
    \end{equation}
the local times of the $i$-th piece. Fix two parameters $\delta,C\in(0,\infty)$
and estimate
    \begin{equation}
    \label{lowerb1}
    E\bigl(e^{-\beta H_n}\1_{\{S_n\geq  b\beta^{\frac 13}n\}}\bigr)
    \geq E\Big(e^{-\beta H_n} \prod_{i=1}^{n/T_\beta}
    \bigl[\1_{\ev_i(\delta,T, \beta)}\1_{\ev^\leq_i(\delta,C,T,\beta)}
    \1_{\{S^{\smallsup{i}}_{T_\beta} \geq  b \beta^{\frac 13}T_{\beta}\}}\bigr]\Bigr),
    \end{equation}
where the events $\ev_i(\delta,T,\beta)$ and $\ev^\leq_i(\delta,T,C,\beta)$ are defined by
    \begin{eqnarray}
    \ev_i(\delta,T,\beta)&=&\Bigl\{S^{\smallsup{i}}
    \subset[-\delta \beta^{-\frac 13},
    S^{\smallsup{i}}_{T_\beta}+\delta\beta^{-\frac 13}]\Bigr\},\\
    \ev^\leq_i(\delta,T,C,\beta)&=&\Bigl\{\max_{x\colon|x| \leq
    \delta\beta^{-\frac 13}} \ell^{\smallsup{i}}(x)
    \leq C\beta^{-\frac 13},\max_{x\colon |x-S^{(i)}_{T_\beta}|
    \leq \delta\beta^{-\frac 13}}
    \ell^{\smallsup{i}}(x) \leq C\beta^{-\frac 13}\Bigr\}.
    \label{E<def}
    \end{eqnarray}

\medskip\noindent
{\bf 2.}\ Next, assume that $\delta<bT/2$ (i.e., $\delta\beta^{-\frac 13}<
b\beta^{\frac 13}T_\beta/2$). Then, on the event $\bigcap_{i=1}^{n/T_\beta}
[\ev_i(\delta,T,\beta)\cap \ev^\leq_i(\delta,T,C,\beta)]$, the following hold:
(a) there are no mutual intersections between the pieces unless they are neighbors
of each other; (b) the $i$-th and the $(i+1)$-st piece have mutual intersections
in an interval of length $2\delta \beta^{-\frac 13}$ centered at $S_{iT_\beta}$
only; (c) in this interval the local times of the $i$-th and the $(i+1)$-st piece
are at most $C\beta^{-\frac 13}$, so that the interaction between them satisfies
    \begin{equation}
    \label{estpieces}
    e^{-2\beta \sum_{x} \ell^{\smallsup{i}}(x+S_{(i-1)T_\beta})
    \ell^{(i+1)}(x+S_{iT_\beta})}\geq e^{-4\delta C^2}.
    \end{equation}
Therefore, using \eqref{Hn'def}, together with \eqref{estpieces} and \eqref{E<def}, yields
that on the event $\bigcap_{i=1}^{n/T_\beta}
[\ev_i(\delta,T,\beta)\cap \ev^\leq_i(\delta,T,C,\beta)]$, we have
    \begin{equation}
    e^{-\beta H_n} = e^{-\beta H_n' -2\beta (\ell_n(0)-1)}
    \geq e^{-2 C\beta^{\frac 23}} e^{-\beta H_n'}\geq
    e^{-2 C\beta^{\frac 23}} e^{-4\delta C^2 n/T_{\beta}} \prod_{i=1}^{n/T_{\beta}}
    e^{-\beta H_{T_{\beta}}'(i)},
    \label{estpieces2}
    \end{equation}
where $H_{T_{\beta}}'(i)$ denotes $H_{T_{\beta}}'$ computed for the $i^{\rm th}$ walk
$S^{\smallsup{i}}$. We substitute \eqref{estpieces2} into \eqref{lowerb1}
and note that, after this is done, the pieces are independent. This reasoning
yields
    \begin{equation}
    E\bigl(e^{-\beta H_n}\1_{\{S_n\geq b\beta^{\frac 13}n\}}\bigr)
    \geq e^{-4\delta  C^2  n/T_{\beta}}E\bigl(e^{-\beta H_{T_\beta}'}
    \1_{\ev_1(\delta,T,\beta)}\1_{\ev^\leq_1(\delta,C,T,\beta)}
\1_{\{S^{\smallsup{1}}_{T_\beta} \geq  b \beta^{\frac
13}T_{\beta}\}}\bigr)^{n/T_\beta}.
    \end{equation}

\medskip\noindent
{\bf 3.}\ Next, take logs, multiply by $\beta^{-{\frac 23}}/n= T_\beta/Tn$ and let
$n\to\infty$, to obtain
    \begin{equation}
    \label{estest}
    \beta^{-{\frac 23}}I^+_{\beta}(b\beta^{\frac 13}
    ;b^*\beta^{\frac 13}\sigma^{\frac 23})
    \leq \frac{4\delta C^2}{T}-\frac{1}{T}\log
    E\bigl(e^{-\beta H_{T_\beta}'}\1_{\ev_1(\delta,T,\beta)}\1_{\ev^\leq_1(\delta,C,T,\beta)}
    \1_{\{S^{\smallsup{1}}_{T_\beta} \geq  b \beta^{\frac 13}T_{\beta}\}}\bigr).
    \end{equation}
Let $\beta\da 0$ and use the weak convergence assertions in
\eqref{betatozero}--\eqref{betatozeroWC}, to obtain
    \begin{equation}
    \label{prelimdeltaC}
    \limsup_{\beta \downarrow 0} \bigl[\beta^{-\frac 23} I^+_{\beta}
    (b\beta^{\frac 13};b^*\beta^{\frac 13}\sigma^{\frac 23})\bigr]
    \leq \frac{4\delta C^2}{T}
    - \frac1T\log\widehat E\bigl(e^{-\widehat H^\sigma_T}
    \1_{\widehat \ev(\delta ,T)}\1_{\widehat \ev^\leq(\delta,C,T)}
    \1_{\{B^\sigma_T \geq bT\}}\bigr),
    \end{equation}
where the events $\widehat \ev(\delta,T)$ and $\widehat \ev^\leq(\delta,C,T)$ are defined
in \eqref{Ehatdef}--\eqref{Fhatdeflq}.

\medskip\noindent
{\bf 4.}\ Finally, observe that $\1_{\{B^\sigma_T \geq bT\}}\geq
\1_{\{B^\sigma_T \approx b' T\}}$ for any $b'>b$ and $T$ sufficiently
large (see below \eqref{rateED}). Pick $\mu=\mu_{b'}$ with $\mu_{b'}$
the maximizer in \eqref{LegTra2}, i.e., $\widehat I_1^{\sigma}(b')=\mu_{b'}
b'-\Lambda^+_{\sigma}(\mu_{b'})$. Since $b\geq b^*\sigma^{\frac 23}$ and $b'>b$,
we know that $\mu_{b'}>0$ (recall \eqref{LegTra}). Therefore we may bound
    \begin{equation}
    \label{deviceC}
    \widehat E\bigl(e^{-\widehat H^\sigma_T}
    \1_{\widehat \ev(\delta ,T)}\1_{\widehat \ev^\leq(\delta,C,T)}
    \1_{\{B^\sigma_T \geq bT\}}\bigr)
    \geq e^{-\mu_{b'} b'T+o(T)}\widehat E\bigl(e^{-\widehat H^\sigma_T}e^{\mu_{b'}B^\sigma_T}
    \1_{\widehat \ev(\delta ,T)}\1_{\widehat \ev^\leq(\delta,C,T)}
    \1_{\{B^\sigma_T \approx b'T\}}\bigr).
    \end{equation}
Insert \eqref{deviceC} into \eqref{prelimdeltaC}, let $T\rightarrow \infty$ and
use Proposition~\ref{prop-ingred}(i) (for the Brownian motion with variance
$\sigma^2$ instead of 1), to arrive at
    \begin{equation}
    \label{ratelim}
    \limsup_{\beta \downarrow 0} [\beta^{-\frac 23}
    I^+_{\beta}(b\beta^{\frac 13};b^*\beta^{\frac 13}\sigma^{\frac 23})]
    \leq \mu_{b'} b'-\Lambda^+_{\sigma}(\mu_{b'})=\widehat I_1^{\sigma}(b').
    \end{equation}
Let $b'\da b$ and use the continuity of $\widehat I^\sigma_1$, to complete the proof of
\eqref{betascallow} for $b\geq b^*\sigma^{\frac 23}$.

\medskip\noindent
{\bf 5.}\ The proof of \eqref{betascallow} for $b^{**}\sigma^{\frac 23}<b\leq
b^*\sigma^{\frac 23}$ is analogous. Indeed, \eqref{piece}--\eqref{prelimdeltaC}
give that
    \begin{equation}
    \label{prelimdeltaC2}
    \limsup_{\beta \downarrow 0} \bigl[\beta^{-\frac 23} I^+_{\beta}
    (b\beta^{\frac 13};b^*\beta^{\frac 13}\sigma^{\frac 23})\bigr]
    \leq \frac{4\delta C^2}{T}
    - \frac1T\log\widehat E\bigl(e^{-\widehat H^\sigma_T}
    \1_{\widehat \ev(\delta,T)}\1_{\widehat \ev^\leq(\delta,C,T)}
    \1_{\{0\leq B^\sigma_T\leq bT\}}\bigr),
    \end{equation}
Complete the proof as in \eqref{deviceC}--\eqref{ratelim}, via
$\1_{\{0\leq B^\sigma_T\leq bT\}} \geq \1_{\{B^\sigma_T\approx b'T\}}$
for any $b'<b$ and $T$ sufficiently large, and $\mu_{b'}<0$ for any $b'<b$.
\qed


\section{Proof of Theorem \ref{ratescal}(ii)}
\label{sec-proof(ii)}

In this section we consider the limit $\sigma\to\infty$. Let $(\Ssig_n)_{n\in\N_0}$
be a random walk satisfying \eqref{sigmares}--\eqref{expmomsigma}.

\subsection{Approximative large deviations.}
\label{sec-approxsigma}

Recall \eqref{rate} and \eqref{rate+}. Our main result in this section shows that
the approximative rate function in \eqref{rate+} scales, as $\sigma\to\infty$, to
the rate function for the Edwards model.

    \medskip
    \begin{prop}
    \label{ratescal3}
    Fix $\beta=\infty$. Then, under \eqref{sigmares}--\eqref{expmomsigma},
    \begin{eqnarray}
    \liminf_{\sigma\to\infty}\sigma^{\frac 23}
    I^-_{\infty}\bigl(b\sigma^{\frac 23};b^*\sigma^{\frac 23}\bigr)
    &\geq&\widehat I_1(b),\qquad b\geq 0,\label{rscal}\\
    \limsup_{\sigma\to\infty}\sigma^{\frac 23}
    I^+_{\infty}\bigl(b \sigma^{\frac 23};b^*\sigma^{\frac 23}\bigr)
    &\leq&\widehat I_1(b),\qquad b>b^{**}.\label{rscallow}
    \end{eqnarray}
    \end{prop}

Proposition~\ref{ratescal3} implies Theorem \ref{ratescal}(ii) and Corollary~\ref{asy*}
in the same way as Proposition~\ref{ratescal2} implies Theorem \ref{ratescal}(i) and
Corollary~\ref{asy} (see Section \ref{sec-appro}). We leave this for the reader to
verify.

In the special case of (\ref{rcond}), subject to \eqref{sigmares}--\eqref{expmomsigma},
we know from Theorem~\ref{thm-DJ} that the rate function $I_\beta$ in \eqref{rate} exists
and so we can infer from Proposition~\ref{ratescal3} that
    \begin{equation}
    \lim_{\sigma\to\infty} \sigma^{\frac 23}I_{\infty}\bigl(b\sigma^{\frac 23}\bigr)
    = \widehat I_1(b), \qquad b>b^{**}.
    \end{equation}
Again, we leave open the convergence for $0\leq b \leq b^{**}$.

\subsection{Proof of Proposition \ref{ratescal3}.}
\label{sec-ratescal2*}

Like in Section~\ref{sec-ratescal2}, we decompose the path into pieces to which an
appropriate weak convergence assertion can be applied, which is in this case
\eqref{rtoinfty}. The arguments are similar and again revolve around controlling
the interaction between neighboring pieces. However, it turns out to be more
difficult to handle the mutual {\em avoidance\/} of neighboring pieces than
to handle their mutual {\em intersection local times\/} as in
Section~\ref{sec-ratescal2}. In order to overcome this problem, we use a technique
that is reminiscent of the so-called ``lace expansion''. Throughout the sequel
we write ``$(S_i)_{i=0}^n$ is SAW'' if $S_i\not=S_j$ for all $0\leq i<j\leq n$.

\subsubsection{Proof of (\ref{rscal})}
\label{upperbound*}

\medskip\noindent
{\bf 1.}
Fix $b\geq b^*$ and recall that
     \begin{equation}
     I^+_{\infty}(b\sigma^{\frac 23};b^*\sigma^\frac 23)
     =-\liminf_{n\to\infty}\frac 1n\log P\bigl((\Ssig_j)_{j=0}^n \mbox{ is SAW},
     |\Ssig_n|\geq b\sigma^\frac 23 n\bigr).
     \end{equation}
Instead of (\ref{momZ}), now consider
    \begin{equation}
    \label{momZ*}
    Z_n^{\infty}(\mu) = E\Bigl(e^{\mu \sigma^{-\frac43} \Ssig_n}
    \1_{\{(\Ssig_j)_{j=0}^n \mbox{{\tiny is SAW}}\}}\Bigr),\qquad n\in\N,\,\mu\in\R.
    \end{equation}
Cut the path into $n/T_\sigma$ pieces of length $T_\sigma=\sigma^{\frac 23}T$.
(To simplify the notation, assume that both $n/T_\sigma$ and $T_\sigma$ are integers.)
For $\mu>0$, we estimate, like in \eqref{esti1}--\eqref{esti1a},
    \begin{equation}
    P\bigl((\Ssig_j)_{j=0}^n \mbox{ is SAW},
    \Ssig_n\geq b\sigma^{\frac 23}n\bigr)
    \leq e^{-\mu b \sigma^{-\frac 23}n}\,Z_n^{\infty}(\mu)
    \leq e^{-\mu b \sigma^{-\frac 23}n}\,
    [Z_{T_\sigma}^{\infty}(\mu)]^{n/T_\sigma}.
    \end{equation}

\medskip\noindent
{\bf 2.}\ The following lemma is the analogue of Lemma~\ref{lem-wcbeta} needed here.

    \medskip
    \begin{lemma}
    \label{lem-wcsigma}
    Assume \eqref{sigmares}--\eqref{expmomsigma}. Then, for any $\mu\in\R$,
        \begin{equation}
        \lim_{\sigma\rightarrow \infty} Z_{T_\sigma}^{\infty}(\mu)
        =\widehat E(e^{-\widehat H_T} e^{\mu B_T}).
        \end{equation}
    \end{lemma}

\begin{proofsect}{Proof of Lemma \ref{lem-wcsigma}.}
As in the proof of Lemma \ref{lem-wcbeta}, it suffices to show that
    \begin{equation}
    \limsup_{\sigma\rightarrow \infty}
    E(e^{2\mu\sigma^{-\frac 43}\Ssig_{T_\sigma}})<\infty.
    \end{equation}
Denote the moment generating function of $\Ssig_1/\sigma$ by
$\varphi_{\sigma}(t)=E(e^{t\Ssig_1/\sigma})$. Then
    \begin{equation}
    \label{lapo1*}
    E(e^{2\mu\sigma^{-\frac 43}\Ssig_{T_\sigma}})
    =\varphi_\sigma(2\mu\sigma^{-\frac13})^{T_{\sigma}}.
    \end{equation}
By \eqref{expmomsigma}, the right-hand side is finite for $\sigma$ large enough.
By (\ref{sigmares})(a) we have, uniformly in $\sigma\geq 1$,
    \begin{equation}
    \label{lapo2*}
    \varphi_\sigma(t)=1+\frac 12 t^2 +{\Ocal}(|t|^3),\qquad t\to 0.
    \end{equation}
Put $t=2\mu\sigma^{-\frac13}$ and combine \eqref{lapo1*}--\eqref{lapo2*}, to get
    \begin{equation}
    E\bigl(e^{2\mu\sigma^{-\frac 43}S_{T_\sigma}}\bigr)
    \leq e^{T_\sigma[\frac 12 t^2 + {\Ocal}(1/\sigma)]}
    =e^{2T\mu^2[1+{\Ocal}(\sigma^{-\frac13})]}, \qquad \sigma\to\infty.
    \end{equation}
\end{proofsect}
\qed

\medskip\noindent
{\bf 3.}\ The details of the remainder of the proof are the same as in
Section~\ref{upperbound}, via Lemma~\ref{lem-wcsigma} instead of
Lemma~\ref{lem-wcbeta}. This completes the proof for $b\geq b^*$.
The proof for $0\leq b \leq b^*$ is analogous.
\qed

\subsubsection{Proof of (\ref{rscallow})}
\label{sec-lbAl}

\medskip\noindent
{\bf 1.}\ Fix $b\geq b^*$. Pick any $b'>b$, fix $\sigma,T>0$, and
put $\gamma^{\smallsup{n}}=\gamma_T\sigma^{\frac 23}n/T$.
Then, for $\mu>0$ and $T$ large enough, we have
        \begin{equation}
        \label{sigma1a}
        \1_{\{S_n\geq b\sigma^{\frac 23}n\}}\geq\1_{\{|S_n-b'\sigma^{\frac 23}n|\leq
        \gamma^{\smallsup{n}}\}}
        e^{\mu\sigma^{-\frac 43}[S_n-b'\sigma^{\frac23}n-\gamma^{\smallsup{n}}]}.
        \end{equation}
This implies the lower bound
        \begin{equation}
        \label{sigma1}
        \begin{aligned}
        \sigma^{\frac 23}&\liminf_{n\to\infty}\frac 1n\log
        P\left((S_j)_{j=0}^n\mbox{ is SAW}, S_n\geq b\sigma^{\frac 23}n\right)\\
        &\geq -\mu b'-\mu\frac{\gamma_T}T+\sigma^{\frac 23}\liminf_{n\to\infty}\frac 1n\log
        E\Big(\1_{\{(S_j)_{j=0}^n\mbox{{\tiny is SAW}}\}}e^{\mu\sigma^{-\frac 43}S_n}
        \1_{\{|S_n-b'\sigma^{\frac 23}n|\leq \gamma^{\smallsup{n}}\}}\Big).
        \end{aligned}
        \end{equation}
To handle the expectation in the right-hand side, we estimate
        \begin{equation}
        \label{sigma1b}
        e^{\mu\sigma^{-\frac 43}S_n}
        \1_{\{|S_n-b'\sigma^{\frac 23}n|\leq
        \gamma^{\smallsup{n}}\}}
        \geq\prod_{i=1}^{n/T_\sigma}\Bigl[e^{\mu\sigma^{-\frac 43}S^{\smallsup{i}}_{T_\sigma}}
        \1_{\{|\sigma^{-\frac43}S^{\smallsup{i}}_{T_\sigma}-b'T|\leq \gamma_T\}}
        \1_{\ev_i(\delta,T,\sigma)}\Bigr],
        \end{equation}
where we use the definition \eqref{piece} of the shifted $i$-th piece with $T_\beta$
replaced by $T_\sigma$, abbreviate $S^{\smallsup{i}}=(S^{\smallsup{i}}_j)_{j=0}^{T_\sigma}$,
and introduce the event
        \begin{equation}
        \ev_i(\delta,T,\sigma) = \Bigl\{S^{\smallsup{i}}
        \subset\bigl[-\delta\sigma^\frac 23,
        S^{\smallsup{i}}_{T_\sigma}+\delta\sigma^\frac 23\bigr]\Bigr\}.
        \end{equation}

\medskip\noindent
{\bf 2.}\ Assume that  $\delta<bT/2$. On the event $\bigcap_{i=1}^{n/T_\sigma}
\ev_i(\delta,T,\sigma)$, the pieces $S^{\smallsup{i}}$, $i=1,\dots,n/T_\sigma$,
have no mutual intersection, unless they are neighbors of each other. Hence,
we only need to estimate the interaction between the neighboring pieces. More
precisely, $(S_j)_{j=0}^n$ is SAW as soon as all the pieces $S^{\smallsup{i}}$
are SAW and neighboring pieces do not overlap in more than their connecting point.
Introduce the indicator $U_i$ of the event that the $i$-th and the $(i+1)$-st
piece intersect each other in more than their connecting point:
        \begin{equation}
        U_i(T,\sigma)  = \begin{cases}
        1  &\mbox{if } (S_j)_{j=(i-1)T_\sigma}^{iT_\sigma}
        \cap (S_j)_{j=iT_\sigma}^{(i+1)T_\sigma}
        \not=\{S_{iT_\sigma}\},\\
        0 &\mbox{otherwise}.
        \end{cases}
        \end{equation}
Then we have
        \begin{equation}
        \label{lowbound1a}
        \1_{\{(S_j)_{j=0}^n\mbox{{\tiny is SAW}}\}}
        \prod_{i=1}^{n/T_\sigma}\1_{\ev_i(\delta,T,\sigma)}
        = \prod_{i=1}^{n/T_\sigma}\Bigl[\1_{\{S^{\smallsup{i}}\mbox{{\tiny is SAW}}\}}
        \1_{\ev_i(\delta,T,\sigma)}\Bigr]
        \prod_{i=1}^{n/T_\sigma-1}(1-U_i(T,\sigma)).
        \end{equation}
Using \eqref{sigma1b} and \eqref{lowbound1a}, we obtain the lower bound
        \begin{equation}
        \label{lowbound1}
        E\Big(\1_{\{(S_j)_{j=0}^n\mbox{{\tiny is SAW}}\}}e^{\mu\sigma^{-\frac 43}S_n}
        \1_{\{|S_n-b'\sigma^{\frac 23}n|\leq \gamma^{\smallsup{n}}\}}\Big)\geq
        c_{n/T_\sigma}(\delta,T,\sigma,b',\mu),
        \end{equation}
where
        \begin{equation}
        \label{cdef}
        c_N=c_N(\delta,T,\sigma,b',\mu)=E\Bigl(\prod_{i=1}^{N-1}(1-U_i(T,\sigma))
        \prod_{i=1}^N X_i\Bigr), \qquad N\in\N,
        \end{equation}
with
        \begin{equation}
        X_i=e^{\mu\sigma^{-\frac 43}S^{\smallsup{i}}_{T_\sigma}}\1_{\ev_i(\delta,T,\sigma)}
        \1_{\{|\sigma^{-\frac43}S^{\smallsup{i}}_{T_\sigma}-b'T|\leq \gamma_T\}}
        \1_{\{S^{\smallsup{i}}\mbox{{\tiny is SAW}}\}}.
        \end{equation}

\medskip\noindent
{\bf 3.}\ Next use an expansion argument that is reminiscent of the ``lace expansion
technique'', namely, expand the product $\prod_{i=1}^{N-1}(1-U_i)$  in \eqref{cdef} as
        \begin{equation}\label{expansion}
        \prod_{i=1}^{N-1} (1-U_i) = \sum_{m=1}^N
        \prod_{i=1}^{m-1} (-U_i) \prod_{i=m+1}^{N-1} (1-U_i),
        \end{equation}
where the empty product is defined to be equal to 1. This expansion
has the advantage that every summand splits into a product of two
separated products. Substitute \eqref{expansion} into \eqref{cdef},
to find that
        \begin{equation}
        c_N = \sum_{m=1}^N (-1)^{m-1}
        E\Big(\Big[\prod_{i=1}^{m-1}U_iX_i\Big]\times X_m\times
        \Big[\prod_{i=m+1}^{N-1}(1-U_i)X_i\Big]\times X_N\Big).
        \end{equation}
Since in the $m$-th summand the term $U_m$ is absent, the two factors between the
two pairs of large square brackets are independent: they depend on the path
$(S_j)_{j=0}^n$ up time $mT_\sigma$, respectively, from time $mT_\sigma$ onwards.
Hence, the $c_N$ satisfy the following renewal relation:
        \begin{equation}
        \label{ceq}
        c_N = c_1c_{N-1}+\sum_{m=2}^N (-1)^{m-1}\pi_m c_{N-m},
        \qquad N\in\N,
        \end{equation}
where
        \begin{equation}
        \label{pidef}
        \pi_m = \pi_m(\delta,T,\sigma,b',\mu)
        = E\Big(\prod_{i=1}^{m-1} U_i\prod_{i=1}^{m}X_i\Big).
        \end{equation}

\medskip\noindent
{\bf 4.}\ Use the Cauchy-Schwarz inequality, to estimate
    \begin{equation}
    \pi_m \leq E\Big(\prod_{\stackrel{i=1}{i~{\rm odd}}}^{m-1} U_i \prod_{i=1}^m X_i\Big)^{1/2}
    E\Big(\prod_{\stackrel{i=1}{i~{\rm even}}}^{m-1} U_i \prod_{i=1}^m X_i\Big)^{1/2}
    = (\pi_2^{m/2})^{1/2} c_1 (\pi_2^{(m-2)/2})^{1/2},
    \qquad m\in\N \mbox{ even},
    \end{equation}
and similarly for $m\in\N$ odd. Hence
        \begin{equation}
        \label{pibd}
        \pi_m \leq \eps^{m-1} c_1^m, \qquad m\in\N,
        \end{equation}
where
        \begin{equation}
        \label{epsdef}
        \eps=\frac{\sqrt{\pi_2}}{c_1}.
        \end{equation}

\medskip\noindent
{\bf 5.}\ The following two lemmas give us control over $\eps$ and $c_N$. From now on,
we choose $\mu=\mu_{b'}$ with $\mu_{b'}$ the maximizer in \eqref{LegTra}, i.e.,
$\widehat I(b')=\mu_{b'} b'-\Lambda^+(\mu_{b'})$, which is possible when $b'>b^{**}$
(recall \eqref{LegTra}).

        \medskip
        \bel
        \label{epssmall}
        Fix $b'>b^{**}$. Then
        \begin{equation}
        \label{epslim}
        \lim_{\delta\da 0}\limsup_{T\to\infty}\limsup_{\sigma\to\infty}
        \eps(\delta,T,\sigma,b',\mu_{b'})=0.
        \end{equation}
        \el

        \medskip
        \bel
        \label{renewal}
        For $\eta>0$ sufficiently small the following is true:
        If $\delta,T,\sigma>0$ are chosen such that $\eps=\eps(\delta,T,\sigma,b',\mu_{b'})
        <\eta$, then there are numbers $C,N_0>0$ (depending on $\eps$ and $\eta$ only)
        such that
            \begin{equation}
            \label{cesti}
            c_N \geq C(1-3\eta)^Nc_1^N,\qquad N>N_0.
            \end{equation}
            \el

\medskip\noindent
{\bf 6.}\ Before giving the proof of Lemmas \ref{epssmall}--\ref{renewal}, we complete
the argument. Pick $\eta\in(0,\frac 14)$ so small that Lemma~\ref{renewal} is satisfied
for this $\eta$. According to Lemma~\ref{epssmall}, we may pick $\delta>0$ so small that,
when $T$ is picked sufficiently large, we have $\eps<\eta$ for any sufficiently large
$\sigma$. Hence we may make use of the estimate in \eqref{cesti} for these $T$ and
$\sigma$.

We use \eqref{lowbound1} and Lemma~\ref{renewal} in \eqref{sigma1}, to obtain
        \begin{equation}
        \label{sigma1est}
        \begin{aligned}
        \sigma^{\frac 23}I^+_{\infty}(b\sigma^{\frac 23};b^*\sigma^{\frac 23})
        &=-\sigma^{\frac 23}\liminf_{n\to\infty}
        \frac 1n\log P\bigl((S_j)_{j=0}^n\mbox{ is SAW}, S_n\geq
        b\sigma^{\frac 23}n\bigr)\\
        &\leq \mu_{b'} b'+\mu_{b'}\frac{\gamma_T}T-\liminf_{n\to\infty}
        \frac{\sigma^{\frac 23}}n\log c_{n/T_\sigma}\\
        &\leq \mu_{b'} b'+\mu_{b'}\frac{\gamma_T}T
        -\liminf_{n\to\infty}\frac {\sigma^{\frac 23}}n
        \log\bigl[C(1-3\eta)^{n/T_\sigma}c_1^{n/T_\sigma} \bigr]\\
        &=\mu_{b'} b'+\mu_{b'}\frac{\gamma_T}{T}-\frac{1}{T}\log(1-3\eta)
        -\frac{1}{T}\log c_1.
        \end{aligned}
        \end{equation}
Return to \eqref{rtoinfty} and recall that $\{H_{T_\sigma}=0\}=\{S^{\smallsup{1}}
\mbox{ is SAW}\}$. From the weak convergence assertion in \eqref{rtoinfty} applied
to \eqref{cdef} for $N=1$, in combination with a statement like in Lemma~\ref{lem-wcsigma},
it follows that
        \begin{equation}
        \label{cconv}
        \lim_{\sigma\to\infty}c_1(\delta,T,\sigma,b',\mu_{b'})=
        \widehat E\bigl(e^{-\widehat H_T}e^{\mu_{b'} B_T}\1_{\widehat \ev(\delta,T)}
        \1_{\{B_T\approx b'T\}}\bigr),
        \end{equation}
where $\widehat \ev(\delta,T)$ is the event defined in \eqref{Ehatdef}. Combining
\eqref{sigma1est}--\eqref{cconv}, we obtain
        \begin{equation}
        \label{finish}
        \limsup_{\sigma\to\infty}[\sigma^{\frac 23}I^+_{\infty}(b\sigma^{\frac 23}
        ;b^*\sigma^{\frac 23})]\leq \mu_{b'} b'+\mu_{b'}\frac{\gamma_T}T
        -\frac 1T\log(1-3\eta)
        -\frac 1T\log \widehat E\bigl(e^{-\widehat H_T}e^{\mu_{b'} B_T}
        \1_{\widehat \ev(\delta,T)}
        \1_{\{B_T\approx b'T\}}\bigr).
        \end{equation}
Now let $T\to \infty$ and use (\ref{refine1}) for $C=\infty$, to see that the
right-hand side of \eqref{finish} tends to $\mu_{b'} b' -\Lambda^+(\mu_{b'})$,
which is equal to $\widehat I_1(b')$. Finally, let $b'\da b$ and use the continuity
of $\widehat I_1$ to finish the proof of \eqref{rscallow}.
\qed

\subsection{Proof of Lemma \ref{epssmall}.}

\medskip\noindent
{\bf 1.} Fix $\delta,T$. Introduce the Brownian event
        \begin{equation}
        \label{Eihatdef}
        \widehat \ev_i(\delta,T)=
        \bigl\{B_{[(i-1)T,iT]}\subset[-\delta+B_{(i-1)T},B_{iT}+\delta]\bigr\},
        \qquad i=1,2,
        \end{equation}
and note that $\widehat \ev_i(\delta,T)$ is identical to $\widehat \ev(\delta,T)$
in \eqref{Ehatdef} for the $i$-th piece. Write $U_1$ as $1-(1-U_1)$ in the definition
of $\pi_2$ in \eqref{pidef}, to obtain from \eqref{epsdef} that
        \begin{equation}
        \begin{aligned}
        \eps^2&=\frac1{c_1^2}\Bigl[E(X_1X_2)-
        E\Bigl(\1_{\{S^{\smallsup{1}},S^{\smallsup{2}}
        \mbox{{\tiny avoid each other}}\}}X_1X_2\Bigr)\Bigr]\\
        &=1-\frac1{c_1^2}E\Bigl(\1_{\{(S_j)_{j=0}^{2T_\sigma}\mbox{{\tiny is
        SAW}}\}}
            X_1X_2\Bigr)\Bigr].
        \end{aligned}
        \end{equation}
Now apply the weak convergence statement in \eqref{rtoinfty} and recall \eqref{cconv},
to obtain, analogously to \eqref{cconv}, that
        \begin{equation}
        \label{epsconv}
        \lim_{\sigma\to\infty}\eps^2
        =1-\frac{\widehat E\bigl(e^{-\widehat H_{2T}} e^{\mu B_{2T}}
        \1_{\widehat \ev_1(\delta,T)\cap \widehat \ev_2(\delta,T)}\1_{\{B_T\approx b'T\}}
        \1_{\{B_{2T}- B_T\approx b'T\}}\bigr)}
        {\widehat E\bigl(e^{-\widehat H_T}e^{\mu B_T}\1_{\widehat \ev_1(\delta,T)}
        \1_{\{B_T\geq 0\}}\bigr)^2}(1+o(1)),
        \end{equation}
where $o(1)$ refers to $T\to\infty$.

\medskip\noindent
{\bf 2.} Denote the intersection local time of the $i$-th piece by
$\widehat H_T^{\smallsup{i}}$. Then \eqref{epsconv} reads
        \begin{equation}
        \label{epsconva}
        \lim_{\sigma\to\infty}\eps^2
        =\frac{\widehat E\bigl(\bigl[e^{-\widehat H_T^{\smallsup{1}}
        -\widehat H_T^{\smallsup{2}}}
        -e^{-\widehat H_{2T}} \bigr]e^{\mu B_{2T}}
        \1_{\widehat \ev_1(\delta,T)\cap \widehat \ev_2(\delta,T)}
        \1_{\{B_T\approx b'T\}}\1_{\{B_{2T}- B_T\approx b'T\}}}
        {\widehat E\bigl(e^{-\widehat H_T}e^{\mu B_T}\1_{\widehat \ev_1(\delta,T)}
        \1_{\{B_T\geq 0\}}\bigr)^2}(1+o(1)).
        \end{equation}
Denote the local time of the $i$-th piece by $L^{\smallsup{i}}(T,\cdot)$.
Then, on the event $\widehat \ev_1(\delta,T)\cap \widehat \ev_2(\delta,T)$,
we have
        \begin{equation}
        \label{H2T}
        \widehat H_{2T}=\widehat H_T^{\smallsup{1}}+\widehat H_T^{\smallsup{2}}
        +2\int_{B_T-\delta}^{B_T+\delta}L^{\smallsup{1}}(T,x)L^{\smallsup{2}}(T,x)\,\d x.
        \end{equation}
Now fix a small $\alpha>0$ and introduce the events
        \begin{eqnarray}
        \widehat \ev_1^{\,\geq,+}(\delta,\alpha,T)&=&\bigl\{\max_{x\in[B_T-\delta,B_T+\delta]}
        L^{\smallsup{1}}(T,x)\geq\alpha\delta^{-1/2}\bigr\},\\
        \widehat \ev_2^{\,\geq,-}(\delta,\alpha,T)&=&\bigl\{\max_{x\in[B_T-\delta,B_T+\delta]}
        L^{\smallsup{2}}(T,x)\geq\alpha\delta^{-1/2}\bigr\}.
        \end{eqnarray}
We estimate the right-hand side of \eqref{epsconv} differently on the event
$\widehat \ev_1^{\,\geq,+}\cup \widehat \ev_2^{\,\geq,-}$ and on its complement.
Namely, on the complement of $\widehat \ev_1^{\,\geq,+}\cup \widehat \ev_2^{\,\geq,-}$
we estimate
        \begin{equation}
        \widehat H_{2T}\leq \widehat H_T^{\smallsup{1}}+\widehat
        H_T^{\smallsup{2}} +4\alpha^2,
        \end{equation}
which implies
        \begin{equation}
        \label{esticompl}
        e^{-\widehat H_T^{\smallsup{1}}-\widehat H_T^{\smallsup{2}}}-
        e^{-\widehat H_{2T}}\leq \bigl[1-e^{-4\alpha^2}\bigr]
        e^{-\widehat H_T^{\smallsup{1}}-\widehat H_T^{\smallsup{2}}},
        \end{equation}
while on the event $\widehat \ev_1^{\,\geq,+}\cup \widehat \ev_2^{\,\geq,-}$ we estimate
$-e^{-\widehat H_{2T}}\leq 0$. By symmetry, $\widehat \ev_1^{\,\geq,+}$ and
$\widehat \ev_2^{\,\geq,-}$  have the same probability. Summarizing, we obtain
        \begin{equation}
        \label{epsesti}
        \lim_{\sigma\to\infty}\eps^2 \leq 1-e^{-4\alpha^2}+2
        \left(\frac{\widehat E\bigl(e^{-\widehat H_T}e^{\mu B_T}
        \1_{\widehat \ev(\delta,T)}\1_{\widehat \ev^{\,\geq}(\delta,\alpha,T)}
        \1_{\{B_T\approx b'T\}}\bigr)} {\widehat E\bigl(e^{-\widehat H_T}
        e^{\mu B_T}\1_{\widehat \ev(\delta,T)}\1_{\{B_T\geq0\}}\bigr)}\right)^2(1+o(1)),
        \end{equation}
where we recall that the events $\widehat \ev(\delta,T)$ and $\widehat
\ev^{\,\geq}(\delta,\alpha,T)$ are defined in \eqref{Ehatdef}, respectively,
\eqref{Fhatdefgq}.

\medskip\noindent
{\bf 3.} Let $T\to\infty$ in \eqref{epsesti} and use
Proposition~\ref{prop-ingred}(i--ii), to obtain
    \begin{equation}
    \limsup_{T\to\infty}\lim_{\sigma\to\infty} \eps^2
    \leq 1-e^{-4\alpha^2} + 2\,\frac{K_2(\delta,\alpha)}{K_1(\delta,\infty)}.
    \end{equation}
Let $\delta\da 0$ and use Proposition~\ref{prop-ingred}(iii), to obtain
        \begin{equation}
        \limsup_{\delta\da 0}\limsup_{T\to\infty}\lim_{\sigma\to\infty}
        \eps^2\leq  1-e^{-4\alpha^2}.
        \end{equation}
Let $\alpha\downarrow 0$, to arrive at the assertion in \eqref{epslim}.
\qed

\subsection{Proof of Lemma \ref{renewal}.}
\label{sec-renewal}

\noindent
{\bf 1.}\ Fix $\eta>0$ and $\eps\in(0,\eta)$. Define $z\in(0,\infty)$ by
        \begin{equation}
        \label{zdef}
        1-z=\sum_{m=2}^\infty (-1)^{m-1}\pi_m\Bigl(\frac z{c_1}\Bigr)^m.
        \end{equation}
Equation \eqref{pibd} implies that, for any $z\in(0,\frac 1\eta)$, the modulus
of the right-hand side is bounded above by $\eps z^2/(1-\eps z)\leq \eta z^2/(1-\eta z)$.
Since this function crosses $1-z$ in $z=\frac 1{1+\eta}$ and since, for sufficiently
small $\eta$, its negative value crosses $1-z$ in $\frac 32+\Ocal(\eta)$, there is
indeed a solution $z$ to \eqref{zdef} in $(0,\frac 32+\Ocal(\eta)]$ as $\eta\downarrow 0$.
Assume that $\eta$ is so small that this solution exists and satisfies the estimate
$z^{-1}\geq 1-3\eta$.

\medskip\noindent
{\bf 2.}\ Abbreviate
    \begin{equation}
    \label{Adef}
    A_N=c_N\Bigl(\frac z{c_1}\Bigr)^N,\qquad N\in\N \qquad (A_0=1).
    \end{equation}
We claim that, if $\eta$ is small enough, then there are numbers $K>0$ and $q\in(0,1)$
(depending on $\eta$ only) such that
    \begin{equation}
    \label{delA}
    |A_N-A_{N-1}|\leq K q^N, \qquad N\in\N.
    \end{equation}
The proof of this claim is given in part 3. Because of \eqref{delA},
$A_\infty=\lim_{N\to\infty}A_N\in(0,\infty)$ exists and, for $N$ sufficiently
large, $A_N\geq \frac 12 A_\infty$, which reads $c_N\geq \frac 12 A_\infty
(\frac{c_1}{z})^{N}$. Recall that $z^{-1}\geq 1-3\eta$, to finish
the proof of the lemma.

\medskip\noindent
{\bf 3.}\ The proof of \eqref{delA} goes via induction on $N$. Pick $q=\sqrt{\eta z}$
and assume that $\eta$ is so small that $q<1$ and
   \begin{equation}
   \label{etacond}
   \frac{\eta z^2}{(1-q)(q-\eta z)}\leq \frac 12.
   \end{equation}
Furthermore, pick $K\geq 1$ so large that $1+Kq/(1-q)\leq K(1-\eta z)/2z$.
Then the claim holds for $N=1$, since $|A_1-A_0|=|z-1|\leq \eta z^2/(1-\eta z)
\leq \frac 12 (q-\eta z) \leq Kq$. Assume now that $N>1$ and that the claim
holds for all positive integers $<N$. From this induction hypothesis it follows
that, for every $m=2,\dots,N$,
    \begin{equation}
    \label{Atel}
    |A_{N-1}-A_{N-m}|\leq \sum_{k=2}^m|A_{N-k+1}-A_{N-k}|\leq Kq^N\sum_{k=2}^mq^{-k+1}
    =\frac{Kq}{1-q}q^{N-m}.
    \end{equation}
Estimate, with the help of \eqref{ceq}, \eqref{pibd}, \eqref{zdef}--\eqref{Adef},
the triangle inequality and \eqref{Atel},
    \begin{equation}
    \label{Acalc}
    \begin{aligned}
    |A_N-A_{N-1}|&= \Bigl| c_{N-1}\Bigl(\frac z{c_1}\Bigr)^{N-1}(z-1)
    +\Bigl(\frac z{c_1}\Bigr)^{N}\sum_{m=2}^N(-1)^{m-1}\pi_mc_{N-m}\Bigr|\\
    &=\Bigl| -A_{N-1}\sum_{m=2}^\infty (-1)^{m-1}\pi_m
    \Bigl(\frac z{c_1}\Bigr)^{m}+
    \sum_{m=2}^N(-1)^{m-1}\pi_m\Bigl(\frac z{c_1}\Bigr)^{m}A_{N-m}\Bigr|\\
    &\leq A_{N-1}\Bigl|\sum_{m=N+1}^\infty (-1)^{m-1}\pi_m\Bigl(\frac z{c_1}\Bigr)^{m}\Bigr|
    +\Bigl|\sum_{m=2}^N(-1)^{m-1}\pi_m\Bigl(\frac z{c_1}\Bigr)^m
    \bigl(A_{N-m}-A_{N-1}\bigr)\Bigr|\\
    &\leq \bigl(1+|A_{N-1}-A_0|\bigr)\sum_{m=N+1}^\infty \eps^{m-1} z^m
    +\sum_{m=2}^N\eps^{m-1}z^m|A_{N-m}-A_{N-1}|\\
    &\leq \Bigl(1+\frac{Kq}{1-q}\Bigr)\frac{\eta^Nz^{N+1}}{1-\eta z}
    +\frac{Kq}{1-q} \frac{q^N}\eta
    \sum_{m=2}^\infty \Bigl(\frac{\eta z}q\Bigr)^m\\
    &=\Bigl(1+\frac{Kq}{1-q}\Bigr)\frac{z(\eta z)^N}{1-\eta z}
    +K q^N\frac{\eta z^2}{(1-q)(q-\eta z)}.
    \end{aligned}
\end{equation}
Now recall that $1+Kq/(1-q)\leq K(1-\eta z)/2z$ and recall the estimate in
\eqref{etacond}. Furthermore, observe that $\eta z\leq \sqrt{\eta z} = q$.
This implies that the right-hand side of \eqref{Acalc} is at most $Kq^N$, which
finishes the proof of the induction step.
\qed


\section{Remaining proofs}
\label{sec-proofs}

In this section we prove the remaining results in Section~\ref{results}:
Theorems~\ref{thm-betan}--\ref{thm-twodim}. All the proofs are minor
adaptations of the proof in Section~\ref{sec-proof(i)}.

\subsection{Proof of Theorem~\ref{thm-betan}.}
\label{proof-betan}

The main result proved in this section is the analogue of Proposition~\ref{ratescal2}
for the case where the strength of self-repellence $\beta$ is coupled to the length
of the polymer $n$:

    \medskip
    \begin{prop}
    \label{thm-betan2}
    Assume \eqref{mv}. If $\beta$ is replaced by $\beta_n$ satisfying $\beta_n \to 0$
    and $\beta_n n^{\frac 32} \to\infty$ as $n\to\infty$,
    then
        \begin{eqnarray}
        -\lim_{n\to\infty}\frac{1}{\beta_n^{\frac 23}n}\log
        E\Bigl(e^{-\beta_n H_n}\1_{\{\Ssig_n\geq b
        \beta_n^{\frac 13}n\}}\Bigr) &=& \widehat I^\sigma_1(b),
        \qquad b\geq b^{*}\sigma^{\frac 23},\\
        -\lim_{n\to\infty}\frac{1}{\beta_n^{\frac 23}n}\log
        E\Bigl(e^{-\beta_n H_n}\1_{\{0\leq\Ssig_n\leq b
        \beta_n^{\frac 13}n\}}\Bigr) &=& \widehat I^\sigma_1(b),\qquad b^{**}
        \sigma^{\frac 23}< b\leq b^{*}\sigma^{\frac 23},\\
        -\liminf_{n\to\infty}\frac{1}{\beta_n^{\frac 23}n}\log
        E\Bigl(e^{-\beta_n H_n}\1_{\{0\leq\Ssig_n\leq b
        \beta_n^{\frac 13}n\}}\Bigr) &\leq& \widehat I^\sigma_1(b),
        \qquad 0\leq b\leq b^{**}\sigma^{\frac 23}.
        \end{eqnarray}
    \end{prop}

\noindent
The proof of Proposition~\ref{thm-betan2} is identical to that of
Proposition~\ref{ratescal2} after we replace the double limit $n\rightarrow\infty$,
$\beta\da 0$ (in this order) by the single limit $n\to\infty$ with
the restrictions $\beta_n\to 0$, $\beta_n n^{\frac 32}\rightarrow\infty$.
The latter implies that $T_{\beta_n} = \beta_n^{-\frac 23}T=o(n)$, and it
is actually only this fact that is needed in the proof. Therefore we can
simply copy the proofs in Sections \ref{upperbound}--\ref{lowerbound*}
to derive Proposition~\ref{thm-betan2}. The reader is asked to check the details.
Proposition~\ref{thm-betan2} in turns implies Theorem~\ref{thm-betan}(i).

A similar result holds when $\sigma$ is coupled to $n$ with the restrictions
$\sigma_n\to\infty$, $\sigma_n n^{-\frac 32}\to 0$. The latter implies that
$T_{\sigma_n}=\sigma_n^{\frac 23}T=o(n)$. The result in turn implies
Theorem~\ref{thm-betan}(ii).
\qed

\subsection{Proof of Theorem~\ref{thm-attr}.}
\label{proof-attr}

Define rate functions $I^+_{\beta,\gamma}$ and $I^-_{\beta,\gamma}$ as in
\eqref{rate+} with $\beta H_n$ replaced by $H_n^{\beta,\gamma}$.
Recall \eqref{betagammarestr}. The main result in this section is the
following.

        \medskip
        \begin{prop}
        \label{thm-attr2}
        Fix $\sigma\in(0,\infty)$. Then, under \eqref{mv},
        \begin{eqnarray}
        \liminf\limits_{\beta,\gamma}
        (\beta-\gamma)^{-\frac 23}
        I^-_{\beta,\gamma}\bigl(b (\beta-\gamma)^{\frac 13}
        ;b^*(\beta-\gamma)^{\frac 13}\sigma^{\frac 23}\bigr)
        &\geq&\widehat I^\sigma_1(b),\qquad b\geq 0,\label{attrlower}\\
        \limsup\limits_{\beta,\gamma}(\beta-\gamma)^{-\frac 23}
        I^+_{\beta,\gamma}\bigl(b(\beta-\gamma)^{\frac 13}
        ;b^*(\beta-\gamma)^{\frac 13}\sigma^{\frac 23}\bigr)
        &\leq&\widehat I^\sigma_1(b),\qquad b>b^{**}\sigma^{\frac 23}.
        \label{attrupper}
        \end{eqnarray}
        \end{prop}

\medskip\noindent
Proposition~\ref{thm-attr2} implies Theorem~\ref{thm-attr}, analogously to the
proof in Section~\ref{sec-appro}. We believe that Proposition~\ref{thm-attr2}
and Theorem~\ref{thm-attr} fail without the restrictions on $\beta,\gamma$ in
(\ref{betagammarestr}).

\subsubsection{Proof of \eqref{attrlower}}

Fix $b\geq b^*\sigma^{\frac 23}$. Fix $T>0$ and put $T_{\beta,\gamma}
=T_{\beta-\gamma}=T(\beta-\gamma)^{-\frac 23}$. (Again, assume for notational
convenience that both $T_{\beta-\gamma}$ and $n/T_{\beta-\gamma}$ are integers.)
First note that the interaction in \eqref{Hamiattr} may be written as
        \begin{equation}
        \label{Hident}
        H_n^{\beta,\gamma}
        = (\beta-\gamma) H_n +\frac \gamma2 G_n,
        \end{equation}
where $H_n$ is the interaction of the Domb-Joyce model in \eqref{Hamilt},
and
\begin{equation}
G_n=\sum_{x\in\Z}[\ell_n(x)-\ell_n(x+1)]^2.
\end{equation}
(Absorb the terms $n+1$ in \eqref{Hamilt} and $\beta(n+1)$ in \eqref{Hamiattr}
into the normalization.)
Define
        \begin{equation}
        Y_n^{\beta,\gamma}(b)
        = E\bigl(e^{-H_n^{\beta,\gamma}}\1_{\{S_n\geq b(\beta-\gamma)^{\frac 13}n\}}\bigr).
        \end{equation}
To get the lower bound, simply estimate $H_n^{\beta,\gamma}\geq (\beta-\gamma) H_n$
in \eqref{Hident}, which implies that $Y_n^{\beta,\gamma}(b)\leq Y_n^{\beta-\gamma,0}(b)$.
Hence
        \begin{equation}
        \limsup_{n\to\infty}\frac 1n\log Y_n^{\beta,\gamma}(b)
        \leq \limsup_{n\to\infty}\frac 1n\log E\bigl(e^{-(\beta-\gamma)H_n}
        \1_{\{S_n\geq b(\beta-\gamma)^{\frac 13}n\}}\bigr).
        \end{equation}
The right-hand side is nothing but the approximative rate function $I_\beta^+$ defined
in \eqref{rate+} with $\beta$ replaced by $\beta-\gamma$. Hence, \eqref{attrlower}
follows from \eqref{betascal}.

\subsubsection{Proof of \eqref{attrupper}}

\noindent
{\bf 1.} Like in Section~\ref{lowerbound*}, we first show that
(recall \eqref{estest})
        \begin{equation}
        \label{attr1}
        \begin{aligned}
        (\beta-\gamma)^{-\frac 23}&I_{\beta,\gamma}^+\bigl(b(\beta-\gamma)^{-\frac 23}
        ;b^*(\beta-\gamma)^{\frac 13}\sigma^{\frac 23}\bigr)\\
        &=-(\beta-\gamma)^{-\frac 23}\liminf_{n\to\infty}
        \frac 1n\log Y_n^{\beta,\gamma}(b)\\
        &\leq \frac{4\delta C^2}{T}\frac\beta{\beta-\gamma}
        -\frac 1T\log E\Bigl(e^{-H_{T_{\beta,\gamma}}^{\beta,\gamma}}
        \1_{\ev(\delta,T,\beta-\gamma)}\1_{\ev^\leq(\delta,T,C,\beta-\gamma)}
        \1_{\{S_{T_{\beta-\gamma}}\geq b (\beta-\gamma)^{\frac13}T_{\beta,\gamma}\}}\Bigr).
        \end{aligned}
        \end{equation}
With the help of \eqref{Hident} for $n=T_{\beta-\gamma}$
and the inequality $e^{-x} \geq 1-x$, we estimate
        \begin{equation}
        \label{Hbetagamma}
        \begin{aligned}
        e^{-H_{T_{\beta,\gamma}}^{\beta,\gamma}}&=e^{-(\beta-\gamma)H_{T_{\beta-\gamma}}}
        e^{-\frac \gamma 2 G_{T_{\beta-\gamma}}}
        \geq  e^{-(\beta-\gamma)H_{T_{\beta-\gamma}}}
        \bigl[1-\textstyle{\frac \gamma2}G_{T_{\beta-\gamma}}\bigr]
        \geq e^{-(\beta-\gamma)H_{T_{\beta-\gamma}}}
        -\textstyle{\frac \gamma2}G_{T_{\beta-\gamma}}.
        \end{aligned}
        \end{equation}
As to the second term on the right-hand side of \eqref{Hbetagamma}, in part 2 we
show that
    \begin{equation}
    \label{momentcond}
    \lim_{\beta,\gamma}\textstyle{\frac \gamma2}E\big(G_{T_{\beta-\gamma}}\big)=0,
    \end{equation}
where $\lim_{\beta,\gamma}$ is the limit in \eqref{betagammarestr}. Hence, applying
$\lim_{\beta,\gamma}$ on the right-hand side of \eqref{attr1}, we see that the
remainder of the proof is now the same as in Section~\ref{upperbound} after
\eqref{estest}. Thus, the proof is finished as soon as \eqref{momentcond} is proved.

\medskip\noindent
{\bf 2.}\ In order to prove (\ref{momentcond}), we compute
        \begin{eqnarray}
        \label{diffbound}
        E\Big(\sum_{x\in\Z}
        [\ell_{n}(x)-\ell_{n}(x+1)]^2\Big)
        &=&\sum_{i,j=0}^n \big[2P(S_i=S_j)-P(S_i+1=S_j)-P(S_i-1=S_j)\big]\nonumber\\
        &=& \sum_{i,j=0}^n\big[2P(S_{|i-j|}=0)-P(S_{|i-j|}=1)-P(S_{|i-j|}=-1)\big]\\
        &=& 2n+\sum_{k=1}^n \sum_{j=0}^{n-k}
        \big[2P(S_{k}=0)-P(S_{k}=1)-P(S_{k}=-1)\big]\nonumber\\
        &=& 2n+\sum_{k=1}^n (n-k+1)\big[2P(S_{k}=0)-P(S_{k}=1)-P(S_{k}=-1)\big]\nonumber.
        \end{eqnarray}
We must show that the right-hand side of \eqref{diffbound} is $\Ocal(n)$, because then
\eqref{momentcond} follows via our assumption that $\gamma(\beta-\gamma)^{-\frac 23}\to 0$.
This is shown in Lemma \ref{lem-LCLT} below.
\qed

        \medskip
        \begin{lemma}
        \label{lem-LCLT}
        As $n\to\infty$,
        \begin{equation}\label{BNdef}
        \sum_{k=1}^n (n-k+1)\big[2P(S_{k}=0)-P(S_{k}=1)-P(S_{k}=-1)\big]= \Ocal(n).
        \end{equation}
        \end{lemma}

\begin{proofsect}{Proof of Lemma \ref{lem-LCLT}.}
Let $\phi(t) = E(e^{\i t S_1})$ denote the characteristic function of $S_1$.
We have
        \begin{equation}
        P(S_k=x)=\frac{1}{2\pi} \int_{-\pi}^{\pi} e^{\i tx} \phi(t)^k
        \,\d t,\qquad x\in\Z, k\in\N.
        \end{equation}
In particular,
        \begin{equation}
        \label{start LCLT}
        2P(S_{k}=0)-P(S_{k}=1)-P(S_{k}=-1) =\frac{1}{\pi} \int_{-\pi}^{\pi}
        [1-\cos t] \phi(t)^k \,\d t.
        \end{equation}
Abbreviate the left-hand side of \eqref{BNdef} by $B_n$. Then \eqref{start LCLT}
says that
        \begin{equation}
        B_n= \frac{1}{\pi} \int_{-\pi}^{\pi}
        \Bigl[[1-\cos t] \sum_{k=1}^n (n+1-k)\phi(t)^k\Bigr] \,\d t.
        \end{equation}
We next use that
        \begin{equation}
        \sum_{k=1}^n (n+1-k)\phi^k=\frac{n\phi}{1-\phi} -\phi \frac{1-\phi^n}{[1-\phi]^2},
        \qquad\mbox{on }\{\phi\neq 1\},
        \end{equation}
to arrive at
        \begin{equation}\label{BNident}
        B_n= n\,\frac{1}{\pi} \int_{-\pi}^{\pi}
        \phi(t)\frac{1-\cos t}{1-\phi(t)}\,\d t
        -\frac{1}{\pi} \int_{-\pi}^{\pi}\Bigl[\phi(t)\frac{1-\cos t}{1-\phi(t)}\,
        \frac{1-\phi^n(t)}{1-\phi(t)}\Bigr] \,\d t.
        \end{equation}
For the first term, we use that $|\phi(t)|\leq 1$, $t\in[-\pi,\pi]$, and that
the map $t\mapsto\frac{1-\cos t}{1-\phi(t)}$ is bounded on $[-\pi,\pi]\setminus\{0\}$,
since the only value where $\phi(t)=1$ is $t=0$. This shows that the first term
is of order $\Ocal(n)$. For the second term, we use that
        \begin{equation}
        \Big|\frac{1-\phi^n(t)}{1-\phi(t)}\Big| =\Big|\sum_{k=0}^{n-1} \phi^k(t)\Big|
        \leq n, \qquad t\in[-\pi,\pi]\setminus\{0\},
        \end{equation}
so that also the second term in \eqref{BNident} is of order $\Ocal(n)$.
\end{proofsect}
\qed

\subsection{Proof of Theorem~\ref{thm-twodim}.}
\label{proof-twodim}

Let $I^+_{L}$ and $I^-_{L}$ denote the two approximative rate functions for the endpoint of
the {\it first\/} coordinate, $S_n$, with the convention $e^{-\infty H_n}=\1_{\{H_n=0\}}$,
i.e.,
\begin{equation}
    I^+_{L}(\theta;\widetilde \theta) =
    \begin{cases}
    -\liminf\limits_{n\to\infty}\frac 1n\log P^{L}\bigl(H_n=0,
    \Ssig_n\geq\theta n\bigr)&\mbox{if }\theta\geq \widetilde\theta,\\
    -\liminf\limits_{n\to\infty}\frac 1n\log P^{L}\bigl(H_n=0,
    0\leq\Ssig_n\leq \theta n\bigr)&\mbox{if }0\leq\theta\leq \widetilde\theta,
    \end{cases}
    \end{equation}
and similarly for $I^-_{L}$ with $\limsup$. The result below identifies
the asymptotics of these rate functions in the limit as $n\to\infty$ followed by
$L\to\infty$, and also when the two limits are coupled in a
certain way:

    \medskip
    \begin{prop}
    \label{thm-twodim2}
    Fix $\sigma\in(0,\infty)$ and assume \eqref{mv}.
    \begin{enumerate}
    \item[(i)]
        Then
        \begin{eqnarray}
        \liminf_{L\to\infty}L^{\frac 23}
        I^-_{L}\bigl(b(4L)^{-\frac 13}
        ;b^*(4L)^{-\frac 13}\sigma^{\frac 23}\bigr)
        &\geq&\widehat I_1^\sigma(b),\quad b\geq 0,\label{lapoA}\\
        \limsup_{L\to\infty}L^{\frac 23}
        I^+_{L}\bigl(b (4L)^{-\frac 13}
        ;b^*(4L)^{-\frac 13}\sigma^{\frac 23}\bigr)
        &\leq&\widehat I_1^\sigma(b),\quad b>b^{**}\sigma^{\frac 23}
        \label{lapoB}.
        \end{eqnarray}
    \item[(ii)]
        If $L$ is replaced by $L_n$ satisfying $L_n\to \infty$ and
        $L_n n^{-\frac 32}\to 0$ as $n\to\infty$, then
        \begin{eqnarray}
        -\lim_{n\to\infty}\frac{1}{(4L_n)^{-\frac 23}n}\log
        P^{L_n}\Bigl(H_n=0,\Ssig_n\geq b
        (4L_n)^{-\frac 13}n\Bigr) &=& \widehat I^\sigma_1(b),
        \quad b\geq b^{*}\sigma^{\frac 23},
        \label{twodimlowera}\\
        -\lim_{n\to\infty}\frac{1}{(4L_n)^{-\frac 23}n}\log
        P^{L_n}\Bigl(H_n=0,0\leq\Ssig_n\leq b
        (4L_n)^{-\frac 13}n\Bigr) &=& \widehat I^\sigma_1(b),
        \quad b^{**}\sigma^{\frac 23}< b\leq b^{*}\sigma^{\frac 23},\\
        -\liminf_{n\to\infty}\frac{1}{(4L_n)^{-\frac 23}n}\log
        P^{L_n}\Bigl(H_n=0,0\leq \Ssig_n\leq b
        (4L_n)^{-\frac 13}n\Bigr) &\leq& \widehat I^\sigma_1(b),
        \quad 0\leq b\leq b^{**}\sigma^{\frac 23}.
        \label{twodimlowerb}
        \end{eqnarray}
    \end{enumerate}
    \end{prop}
Analogously to before, Theorem~\ref{thm-twodim} is implied by Proposition~\ref{thm-twodim2}.

\begin{proofsect}{Proof of Proposition~\ref{thm-twodim2}.}

\noindent
{\bf 1.}\ Let us compute the conditional probability of the event $\{H_n=0\}$,
i.e., the path $(X_0,\dots,X_n)$ is self-avoiding, given the path
$\Ssig=(\Ssig_0,\dots,\Ssig_n)$ of the first coordinate. Given $\Ssig$, the
event $\{H_n=0\}$ is equal to the event that $U^L_i\not= U^L_j$ for all
time pairs $0\leq i<j\leq n$ at which $\Ssig_i=\Ssig_j$. Let us denote
by $\ell_n(x)$, $x\in\Z$, the local times of $\Ssig$ as in \eqref{deflt},
and by $i_1^x,\dots,i_{\ell_n(x)}^x$ the times at which $\Ssig$ hits $x$.
Then $\{H_n=0\}$ is the event that, for all $x\in\Z$, the random variables
$U^L_{i^x_1},\dots,U^L_{i^x_{\ell_n(x)}}$ are distinct. Since $U^L_0,\dots,U^L_n$
are i.i.d.\ and uniform on $\{-L,\dots,L\}$, the conditional probability of
this event is easily computed:
        \begin{equation}
        \label{PLwrite}
        P^{L}\bigl(H_n=0\,\big|\,\Ssig\bigr)=\prod_{x\in\Z}\prod_{k=0}^{\ell_n(x)-1}
        \Bigl(1-\frac k{2L+1}\Bigr)=\exp\Bigl\{\sum_{x\in\Z}\sum_{k=0}^{\ell_n(x)-1}
        \log\Bigl(1-\frac k{2L+1}\Bigr)\Bigr\}.
        \end{equation}

\medskip\noindent
{\bf 2.}\ Fix $b\geq b^*\sigma^{\frac 23}$. To prove \eqref{lapoB}, use the
inequality $\log(1-x)\leq -x$ and the fact that $\sum_{k=0}^{l-1}k=\frac 12 l(l-1)$,
to estimate
        \begin{equation}
        \label{compute}
        \begin{aligned}
        P^{L}\bigl(H_n=0,\Ssig_n \geq b(4L)^{-\frac 13}n\bigr)
        &=E^{L}\Bigl( \1_{\{\Ssig_n\geq b(4L)^{-\frac 13}n\}}
        P^{L}\bigl(H_n=0\,\big|\,\Ssig\bigr)\Bigr)\\
        &\leq E^{L}\Bigl( \exp \Bigl\{-\frac 1{4L+2}
        \sum\limits_{x\in\Z}\ell_n(x)[\ell_n(x)-1]\Bigr\}
        \1_{\{\Ssig_n\geq b(4L)^{-\frac 13}n\}}\Bigr)\\
        &= E^{L}\Bigl( e^{-\frac 1{4L+2}H_n}
        \1_{\{\Ssig_n\geq b(4 L)^{-\frac 13}n\}}\Bigr),
        \end{aligned}
        \end{equation}
with $H_n$ denoting the self-intersection local time of $\Ssig$ as in \eqref{Hamilt}.
The right-hand side of \eqref{compute} is nothing but the quantity appearing
in \eqref{rate+} for the Domb-Joyce model with strength of self-repellence
$\beta=\frac 1{4L+2}$. For $0\leq b\leq b^{**}\sigma^{\frac 23}$, the same argument
works with $\leq$ replacing $\geq$. Hence, \eqref{lapoA} directly follows from
Proposition~\ref{ratescal2}.

\medskip\noindent
{\bf 3.}\ Fix $b>b^{**}\sigma^{\frac 23}$. To prove \eqref{lapoB}, we insert the
condition that $\max_{x\in\Z}\ell_n(x)\leq \sqrt{L}$. We then have that, for all
$0\leq k(<\ell_n(x))\leq \sqrt{L}$ and $L$ sufficiently large,
        \begin{equation}
        \log\Bigl(1-\frac k{2L+1}\Bigr) \geq  -\frac k{2L+1}\left(1-\frac
        {k}{L}\right)\geq -\frac k{2L+1}\left(1-\frac{1}{\sqrt{L}}\right),
        \end{equation}
and substituting this into (\ref{PLwrite}) we get that
        \begin{equation}
        P^{L}\bigl(H_n=0,S_n\geq b(4L)^{-\frac 13}n\bigr)
        \geq E^{L}\Bigl( e^{-\frac 1{4L+2}(1-\frac{1}{\sqrt{L}}) H_n}
        \1_{\{S_n\geq b(4 L)^{-\frac 13}n\}}
        \1_{\{\max_{x\in\Z} \ell_n(x) \leq \sqrt{L}\}}\Bigr).
        \end{equation}
Now we can follow the same argument as in Section~\ref{lowerbound*}, noting
that the condition $\max_{x\in\Z} \ell_n(x) \leq \sqrt{L}$ is asymptotically
negligible as $L\rightarrow \infty$.

\medskip\noindent
{\bf 4.}\ The proof for $L=L_n$ is identical to the above proof and relies
on Proposition~\ref{thm-betan2}.
\end{proofsect}
\qed

\section{Discussion.}
\label{sec-disc}

The weak interaction limit results in Section~\ref{mainres}--\ref{twodim} were proved
in Sections~\ref{sec-proof(i)}--\ref{sec-proofs} with the help of a new and
flexible method. The idea was to cut the path into pieces of an appropriately
scaled length, to control the interaction between the different pieces,
and to apply the invariance principle to the single pieces. This method
allowed us to prove scaling of the large deviation rate function for the
empirical drift of the path, which in turn implied the weak interaction limit
results in Section~\ref{mainres}--\ref{twodim}.

Our method has a number of advantages over the approach that was followed in our
earlier work, which relied on a variational representation for the quantities in
the central limit theorem and a functional analytic proof that this variational
representation scales to a limit. Our new method is simple, works for a very large
class of random walks in a variety of self-repelling and self-attracting situations,
and allows for a coupled limit in which $n\to\infty$ and $\beta\da 0$, respectively,
$\sigma\to\infty$ together. We expect that it can be applied to other polymer models
as well, such as branched polymers and heteropolymers, which we hope to investigate
in the future.

The results in Section~\ref{mainres}--\ref{twodim} show universality, in the sense
that the scaling limits do not depend on the details of the underlying random walk
other than its step variance and are all given in terms of the Edwards model.

Two items remained open. First, we did not prove the scaling of the rate function
in the linear regime (recall the remark at the end of Sections~\ref{sec-appro}
and \ref{sec-approxsigma}). In this regime we only derived the upper bound. We
have no doubt that the lower bound can be derived too, but this would require
some further refinements. In particular, in the linear regime the path makes
an overshoot, and we would need to control the interaction between overlapping
pieces in the overshoot. Second, we did not prove the scaling of the variance
in the central limit theorem (recall \eqref{sigscal}). This would require
control of the second derivative of the rate function in its minimum (compare
Theorem~\ref{thm-DJ}(iii) with Theorem~\ref{thm-LDEM}(iii)). We only have good
control over the first derivative of the rate function. The LDP does not imply
the CLT, so even if we had obtained the scaling of the variance, we would not
be able to deduce the CLT anyway.

\end{document}